\documentclass[12pt,reqno]{amsart}
\usepackage{amssymb}
\usepackage{amsmath}
\usepackage{mathabx}
\usepackage{mathtools}
\usepackage{amsthm}
\usepackage{amstext}
\usepackage{amsopn}
\usepackage{mathrsfs}
\usepackage{latexsym}
\usepackage{enumerate}
\usepackage{bm}
\usepackage{cleveref}
\DeclareMathAlphabet{\mathpzc}{OT1}{pzc}{m}{it}
\allowdisplaybreaks
\newtheorem{Definition}{Definition}[section]
\newtheorem{Proposition}{Proposition}[section]
\newtheorem{Lemma}{Lemma}[section]
\newtheorem{Theorem}{Theorem}[section]
\newtheorem{Corollary}{Corollary}[section]
\newtheorem{Remark}{Remark}[section]
\newtheorem{Example}{Example}[section]

\numberwithin{equation}{section}
\begin{document}
%%%%%%%%%%%%%%%%%%%%%%%%%%%%%%%%%%%%%%%%%%%%%%%%%%%%%%%%%%%%
%\bibliographystyle{plain}
%%%%%%%%%%%%%%%%%%%%%%%%%%%%%%%%%%%%%%%%%%%%%%%%%%%%%%%%%%%%%%%%%
\title[Infinitesimal moments in free and c-free probability]
{Infinitesimal moments in free and c-free \\[3pt] probability and Motzkin paths}
\author[R. Lenczewski]{Romuald Lenczewski}
\address{Romuald Lenczewski \newline
Katedra Matematyki, Wydzia\l{} Matematyki, Politechnika Wroc\l{}awska, \newline
Wybrze\.{z}e Wyspia\'{n}skiego 27, 50-370 Wroc{\l}aw, Poland}
\email{Romuald.Lenczewski@pwr.edu.pl}
%%%%%%%%%%%%%%%%%%%%%%%%%%%%%%%%%%%%%%%%%%%%%%%%%%%%%%%%%%%%%%%%
\begin{abstract}
Infinitesimal moments associated with infinitesimal freeness and 
infinitesimal conditional freeness are studied.
For free random variables, we consider continuous deformations 
of moment functionals associated with Motzkin paths $w\in \mathpzc{M}$, which provide 
a decomposition of their moments, and we compute their derivatives at zero. We show that the
first-order derivative of each functional vanishes unless the path has exactly one local maximum. 
Geometrically, this means that $w$ is a pyramid path, which is consistent with the characteristic formula 
for alternating moments of infinitesimally free centered random variables.
In this framework, infinitesimal Boolean independence is also obtained and it corresponds to flat paths.
A similar approach is developed for infinitesimal conditional freeness, for which we show that the only moment functionals 
that have a non-zero first-order derivative are associated with concatenations of a pyramid path and a flat path.
This charaterization leads to a Leibniz-type definition of infinitesimal conditional freeness at the level of moments.

\end{abstract}
\maketitle

%\tableofcontents

\section{Introduction}

Free probability, originated by Voiculescu, laid the foundations of modern 
noncommutative probability theory based on the notion of freeness, or free independence \cite{V1}. 
Of particular importance is his discovery that independent random matrices are asymptotically free \cite{V2}. 

Infinitesimal freeness, introduced by Belinschi and Shlyakhtenko \cite{BSh}, plays an analogous role at the next order, as it 
characterizes the first-order corrections to asymptotic freeness arising in random matrix models \cite{Sh}.
Freeness of type B, introduced by Biane, Goodman, and Nica \cite{BGN}, provided a combinatorial framework that later turned out to be closely related to infinitesimal freeness as shown by F\'evrier and Nica \cite{FN}. 
For recent random matrix models in which both infinitesimal freeness and conditional freeness arise naturally,
see C\'ebron, Dahlqvist and Gabriel \cite{CDG}.

The concept of infinitesimal freeness of a family of unital subalgebras $({\mathcal A}_{i})_{i\in I}$ 
of a given unital algebra ${\mathcal A}$ is based on replacing a single normalized linear functional  
$\varphi:{\mathcal A}\rightarrow {\mathbb C}$, with respect to which this family is free,
by a family of deformed normalized linear functionals $(\varphi_{t})_{t\geq 0}$. To study infinitesimal freeness,
it suffices to take deformed functionals $\varphi_{t}=\varphi+t\varphi'+o(t)$ and to compute the derivatives of moments at $t=0$.
The resulting Leibniz-type formula for infinitesimal moments of variables that are free with respect to
$\varphi$ is of the form
\begin{equation}\label{eq:1.1}
\varphi'(a_1\cdots a_n)=\sum_{k=1}^{n}\varphi'(a_k)\varphi(a_1\cdots a_{k-1}a_{k+1}\cdots a_n),
\end{equation}
where $\varphi'$ is another linear functional on ${\mathcal A}$ satisfying $\varphi'(1)=0$,
and where $a_1\in {\mathcal A}_{i_1}, \ldots, a_n\in {\mathcal A}_{i_n}$ are $\varphi$-centered 
and the indices are alternating, that is $i_1\neq \cdots \neq i_n$. 

This formula is equivalent to the characteristic formula for infinitesimal alternating moments of random variables
used by F\'evrier and Nica \cite{FN} (a related type~B version appeared earlier in \cite{BGN}), namely
\begin{equation}\label{eq:1.2}
\varphi'(a_1\cdots a_n)=\varphi(a_1a_n)\cdots \varphi(a_{m-1}a_{m+1})\varphi'(a_m)
\end{equation}
whenever $n=2m-1$ and $i_1=i_n, \ldots, i_{m-1}=i_{m+1}$, while the remaining alternating moments of 
$\varphi$-centered variables with respect to $\varphi'$ vanish. The appearance of such an expression is 
intriguing and suggests a deeper combinatorial structure. Note in this context that the above product of moments 
is equal to the product of marginal moments with respect 
to the normalized linear functionals $\varphi_i=\varphi|_{{\mathcal A}_{i}}$.

In this work, we show that the special form of the above formula, in which
subalgebra labels $i_1, \ldots, i_n$ have to match pairwise except for 
the middle singleton, is linked to pyramid Motzkin paths.
The decomposition of moments of free random variables in terms of moments associated 
with Motzkin paths were subject of our investigations in \cite{L4,L5}.
In particular, we introduced multilinear moment functionals
\begin{equation}\label{eq:1.3}
\Phi_w(a_1, \ldots, a_n)=\Phi(a_1(j_1)\cdots a_n(j_n))
\end{equation}
where $w=j_1\cdots j_n$ is a Motzkin word, the variables $a_k(j_k)$ are copies of 
the $a_{k}\in {\mathcal A}_{i_k}$ called orthogonal replicas and $\Phi$ is a suitably 
defined tensor product functional. These functionals, called {\it Motzkin functionals} in \cite{L4},
contain detailed information about the moments of free random variables 
(we denoted them $\psi(w)$).
We also proved in \cite{L5} that free cumulants 
can be decomposed in terms of the associated Motzkin cumulants.
The framework of Motzkin functionals 
turns out to be also very useful in the infinitesimal context.

Namely, if $({\mathcal A}_{i}, \varphi_i)_{i\in I}$ is a family of 
noncommutative probability spaces and $\star_{i\in I}{\mathcal A}_{i}$ 
is their free product with identified units \cite{Av, V1}, then the linear functional $\varphi$ on the free 
product without identification of units $*_{i\in I}{\mathcal A}_{i}$ given by 
\setlength{\belowdisplayshortskip}{\baselineskip}
\begin{equation}\label{eq:1.4}
\varphi (a_1\cdots a_n)=\sum_{w\in\mathpzc{M}_{n}}\Phi_w(a_1,\ldots, a_n)
\end{equation}
satisfies the equation $\varphi=(\star_{i\in I} \varphi_i)\circ \tau$,
where $\tau$ is the unit identification map. Here,
the variables $a_1\in \mathcal{A}_{i_1}, \ldots, a_n\in \mathcal{A}_{i_n}$ are assumed to 
have alternating labels.
 
We can now differentiate the right-hand side of the above formula. 
For that purpose, we replace each $\Phi_w$ by its continuous deformation
$\Phi_{w,t}$, which amounts to replacing the family of functionals 
$(\varphi_i)_{i\in I}$ from which we construct $\Phi$ by
its continuous deformation $(\varphi_{i,t})_{i\in I}$. 
Then we evaluate the derivative
{
\setlength{\belowdisplayshortskip}{\baselineskip}
\begin{equation}\label{eq:1.5}
\Phi_{w}'(a_1,\ldots, a_n):=
\left.\left[\frac{d}{dt}\Phi_{w,t}(a_1,\ldots ,a_n)\right]\right|_{t=0}
\end{equation}
}%
for all $w\in \mathpzc{M}_n$ and $n\in {\mathbb N}$. This
allows us to find a decomposition formula for  
$\varphi'$ analogous to \eqref{eq:1.4}. 
It is worth to remark that
the tensor product of functionals admits a 
canonical differentiation.

Moreover, the derivatives of moments defined by \eqref{eq:1.5} have strong `selection rules'.
In particular, if the arguments are alternating $\varphi$-centered variables, 
then all `infinitesimal Motzkin functionals' $\Phi'_{w}$ vanish except for the ones corresponding 
to the {\it pyramid Motzkin paths} identified with 
reduced Motzkin words of the form 
\[
w=j_1\cdots j_n=1\cdots (m-1)m(m-1)\cdots 1
\]
where $n=2m-1$, and these assume the form
\[
\Phi_{w}'(a_{1},\ldots, a_{n})=\varphi_{i_1}(a_1a_n)\cdots \varphi_{i_{m-1}}(a_{m-1}a_{m+1})\varphi'_{i_m}(a_m)
\]
whenever $i_1=i_n, \ldots, i_{m-1}=i_{m+1}$. This is the expression that is equivalent to
the non-Leibniz definition of infinitesimal freeness \eqref{eq:1.2}.
Therefore, it is natural to expect that the concept of infinitesimal freeness is related to the Motzkin 
path decomposition of the free product of functionals and of its derivative.

The nonvanishing of $\Phi_{w}'(a_1,\ldots, a_n)$ 
is related to the analytic properties  
of the Motzkin paths in the sense that local maxima of a given path
and their number determine whether $\Phi_{w}'(a_1,\ldots, a_n)$ vanishes or not. 
We should add that a local maximum is understood in the weak sense, namely if $w=j_1\cdots j_n$ and
\[
j_{k-1}\leq j_k\geq j_{k+1}
\]
then the Motzkin word $w$ (and thus the corresponding Motzkin path) has a local maximum at $(k,j_k)$ (if $k=1$ or $k=n$, then 
only one inequality is required). The pyramid Motzkin paths are very special in the class of Motzkin paths: 
they are the only Motzkin paths that have one local maximum and this is the reason why the associated infinitesimal Motzkin 
functionals $\Phi_{w}'$ do not vanish. 

A similar phenomenon occurs when computing higher order derivatives, which allows us 
to study higher order infinitesimal freeness \cite{F} from the perspective of Motzkin paths
and their geometric (or, analytic) properties. For instance, the $m$-th order 
derivative of $\Phi_{w,t}(a_1,\ldots ,a_n)$ at $t=0$ vanishes if $w$
has more than $m$ local maxima.

The infinitesimal Boolean independence is included in the above framework in a very simple fashion
since it corresponds to flat Motzkin paths corresponding to words $w=1^{n}$. 
In this case, we get 
\begin{equation}\label{eq:1.6}
\Phi_{1^{n}}'(a_1, \ldots, a_n)=\sum_{k=1}^{n}\varphi_{i_1}(a_1)\cdots \varphi_{i_k}'(a_k)\cdots 
\varphi_{i_n}(a_n)
\end{equation}
for alternating arbitrary variables $a_1\in {\mathcal A}_{i_1}, \ldots, a_n\in {\mathcal A}_{n}$.
This is the Boolean analog of the Leibniz-type formula \eqref{eq:1.1}.
Let us also remark that infinitesimal monotone independence for two algebras can also be obtained within the above framework; 
however, we do not treat it in this paper, as extending the result to the general case would likely require further generalization.

Moreover, we also show that this approach can be generalized
to the framework of infinitesimal conditional freeness.
This is an infinitesimal version of conditional freeness (c-freeness), the latter originally
introduced and studied by Bożejko, Leinert and Speicher 
\cite{BLS}, where a unital algebra is equipped with a pair of 
normalized linear functionals $(\varphi, \psi)$. We will adopt the definition of conditional freeness of the form
similar to freeness, namely
\begin{equation}\label{eq:1.7}
\varphi(a_{1}\cdots a_n)=0
\end{equation}
whenever $a_1\in {\mathcal A}_{i_1}, \ldots, a_n\in {\mathcal A}_{i_n}$, where 
$i_1\neq \ldots \neq i_n$ and $a_1$ is $\varphi$-centered, whereas $a_k$ is $\psi$-centered 
for $k=2, \ldots, n$. This definition is convenient since it allows us to 
use Motzkin-type functionals to formulate the moment-based definition of infinitesimal conditional 
freeness. 

In the infinitesimal c-free case, we take two families of deformed functionals $(\varphi_{i,t})$ and $(\psi_{i,t})$ for each $i\in I$, and
we suitably modify the deformed Motzkin functionals and their derivatives.
We show that the class of Motzkin paths that give contributions to the first derivative 
of the c-free product of functionals consists of concatenations $w_1w_2$, where $w_1$ is 
a pyramid path and $w_2$ is a flat path, including the cases when one of them is an empty path.
This allows us to formulate a Leibniz-type definition for $\varphi'$-moments 
for infinitesimal conditional freeness with respect to the quadruple of functionals
$(\varphi, \varphi', \psi, \psi')$. 
Namely, for any 
$a_1\in {\mathcal A}_{i_1}^{\circ}$, $a_2\in {\mathcal A}_{i_2}^{\square}, \ldots, a_n\in {\mathcal A}_{i_n}^{\square}$,
we have
\begin{equation}\label{eq:1.8}
\varphi'(a_1\cdots a_n)=\varphi'(a_1)\varphi(a_2 \cdots a_n)+
\sum_{m=2}^{n}\psi'(a_{m})\varphi(a_1\cdots a_{m-1}a_{m+1}\cdots a_n)
\end{equation}
whenever $i_1\neq \cdots \neq i_{n}$, where ${\mathcal A}_{i}^{\circ}={\mathcal A}_{i}\cap {\rm Ker}(\varphi)$
and ${\mathcal A}_{i}^{\square}={\mathcal A}_{i}\cap {\rm Ker}(\psi)$. 
This is a definition of infinitesimal conditional freeness formulated 
directly at the level of moments and a c-free analog of \eqref{eq:1.1}.
The associated cumulants will be studied in a forthcoming paper.

Let us remark that infinitesimal versions of noncommutative independence have so far been studied primarily by means of cumulants \cite{BSh, BGN, FMNS, FN, F, FH} or within an operator-valued framework \cite{T}. In this paper, we develop an alternative approach based on moment computations. This perspective yields new formulations and results that do not arise naturally in cumulant-based treatments. Our method is based on tensor product constructions originating in \cite{L1, L2} and further developed in \cite{L4, L5}. Although these constructions were not originally formulated in infinitesimal terms, we show that they inherently encode infinitesimal information when viewed in light of the work of Shlyakhtenko \cite{Sh} and Mingo and Tseng \cite{MT}, since the idempotents appearing in our constructions play a role analogous to their {\it infinitesimal idempotents}. 
A key conceptual contribution of this paper is the systematic introduction of infinitesimality at two distinct levels: the classical level, through continuous deformations of functionals, and the noncommutative level, through families of infinitesimal idempotents.

The paper is organized as follows. In Section 2 we discuss the concept of infinitesimal freeness. The combinatorics of Motzkin paths 
based on level return blocks and singletons is described in Section 3. Infinitesimal deformations of Motzkin functionals are introduced 
in Section 4. Computing moments of deformed Motzkin functionals is carried out in Section 5 and Section 6.
The main results of the paper provide a decomposition of infinitesimal moments 
in terms of infinitesimal Motzkin functionals and are contained in Section 7 
(infinitesimal freeness), Section 8 (higher order infinitesimal freeness) and Section 9 (infinitesimal conditional freeness). 

\section{Infinitesimal freeness and infinitesimal idempotents}

By a {\it noncommutative probability space} we understand the pair $({\mathcal A}, \varphi)$, where 
${\mathcal A}$ is a unital algebra over ${\mathbb C}$ and $\varphi:{\mathcal A}\rightarrow {\mathbb C}$ 
is a linear functional with $\varphi(1_{\mathcal A})=1$. 
In turn, by an {\it infinitesimal noncommutative probability space} we understand the triple $({\mathcal A}, \varphi, \varphi')$, where
$({\mathcal A}, \varphi)$ is a noncommutative probability space and $\varphi':{\mathcal A}\rightarrow {\mathbb C}$ is a linear 
functional such that $\varphi'(1)=0$. In both cases, an element $a\in {\mathcal A}$ will be called {\it centered} 
or {\it $\varphi$-centered} if $\varphi(a)=0$.

Let us give a Leibniz-type definition of infinitesimal freeness at the level of moments \cite{FN}.

\begin{Definition}
{\rm Let $({\mathcal A}, \varphi, \varphi')$ be an infinitesimal noncommutative probability space and 
let $\{{\mathcal A}_{i}: i\in I\}$ be a family of unital subalgebras of ${\mathcal A}$. This family 
will be called {\it infinitesimally free} with respect to $(\varphi, \varphi')$ if two conditions are satisfied:
\begin{enumerate}
\item
it is free with respect to $\varphi$, namely
\begin{equation}\label{eq:2.1}
\varphi(a_1\cdots a_n)=0
\end{equation}
whenever $a_1\in {\mathcal A}_{i_1}, \ldots, a_n\in {\mathcal A}_{i_n}$ are centered and $i_1\neq \cdots \neq i_n$,
\item
it holds that
\begin{equation}\label{eq:2.2}
\varphi'(a_1\cdots a_n)=\sum_{k=1}^{n}\varphi'(a_k)\varphi(a_1\cdots a_{k-1}a_{k+1}\cdots a_n)
\end{equation}
whenever $a_1 \in {\mathcal A}_{i_1}, \ldots, a_n\in {\mathcal A}_{i_n}$ are centered and $i_1\neq \cdots \neq i_n$.
\end{enumerate}
We say that the family of subsets $\{J_{i}: i\in I\}$ of ${\mathcal A}$ is infinitesimally free with respect to $(\varphi, \varphi')$ 
if the family of unital algebras generated by these subsets $\{{\rm alg}(J_{i}):i\in I\}$
is infinitesimally free with respect to $(\varphi, \varphi')$.}
\end{Definition}

The above definition is equivalent to that giving the explicit formula for mixed moments. 
We state it below in the form of a proposition.

\begin{Proposition}
{\rm Let $({\mathcal A}, \varphi, \varphi')$ be an infinitesimal noncommutative probability space and 
let $\{{\mathcal A}_{i}:\,i\in I\}$ be a family of unital subalgebras of ${\mathcal A}$. This family 
is {\it infinitesimally free} with respect to $(\varphi, \varphi')$ if and only if 
\begin{enumerate}
\item
it is free with respect to $\varphi$,
\item
the mixed alternating moments of centered elements vanish unless $n=2m-1$ for some $m$ and 
$i_1=i_n, \ldots, i_{m-1}=i_{m+1}$, in which case 
\begin{equation}\label{eq:2.3}
\varphi'(a_1\cdots a_n)=
\varphi(a_1a_n)\cdots \varphi(a_{m-1}a_{m+1})\varphi'(a_m)
\end{equation}
whenever $a_1 \in {\mathcal A}_{i_1}, \ldots, a_n\in {\mathcal A}_{i_n}$ are centered and $i_1\neq \cdots \neq i_n$.
\end{enumerate}}
\end{Proposition}

A basic special case of infinitesimal freeness was studied by Mingo and Tseng \cite{MT}, who
considered two subalgebras, with one of them generated by a special idempotent operator $j$. 
Thus, let $({\mathcal A}, \varphi, \varphi')$ be an infinitesimal noncommutative probability 
space, with a unital subalgebra ${\mathcal B}$ and 
let $j\in {\mathcal A}$ be an idempotent such that $\varphi(j)=0$ and $\varphi'(j)\neq 0$. 
This element $j$ is called an {\it infinitesimal idempotent} 
(a slightly more general $j$ such that $\varphi(j^{n})=0$ for all $n\geq 1$ is called an {\it infinitesimal operator}).
Without loss of generality, we will assume that $\varphi'(j)=1$ since
if $\varphi'(j)\neq 1$ we can take $j/\varphi'(j)$. 

If we assume that the family $\{{\mathcal B}, {\rm alg}(\{j\})\}$ is infinitesimally free with respect to $(\varphi, \varphi')$,
then the alternating mixed moments of elements of ${\mathcal B}$ and an infinitesimal idempotent
$j$ under $\varphi'$ assume a nice form. The result given below is not as general as that in \cite{MT}
(see Proposition 5.3), which contains combinatorial formulas for moments involving both $j$ and $j^{\perp}$.
The authors also assumed that $a_1, \ldots, a_n$ are not centered, but there is no loss of generality 
in considering the centered case. For the readers' convenience we also provide a simple proof of our 
version given below. 

\begin{Proposition}
Let $({\mathcal A}, \varphi, \varphi')$ be an infinitesimal noncommutative probability space, with a unital subalgebra 
${\mathcal B}$ and an infinitesimal idempotent $j\in {\mathcal A}$ with $\varphi'(j)=1$. If $\{{\mathcal B}, {\rm alg}(\{j\})\}$ is infinitesimally 
free with respect to $(\varphi, \varphi')$, then it holds that
\begin{enumerate}
\item
for centered $a_1, \ldots, a_n\in{\mathcal B}$, where $n\in {\mathbb N}$,
\[
\varphi'(j^{\alpha}a_1j\cdots ja_{n}j^{\beta})=0
\]
whenever $\alpha, \beta\in \{0,1\}$ unless $\alpha=\beta=0$ and $n=2$,
\item
for centered $a_1, a_2\in \mathcal{B}$,
\[
\varphi'(a_1ja_2)=\varphi(a_1a_2).
\]
\end{enumerate}
\end{Proposition}
{\it Proof.}
If there is exactly one $j$ under $\varphi'$, we either have 
\[
\varphi'(ja_1)=\varphi'(a_1j)=\varphi'(j)\varphi(a_1)=0,
\] 
since $a_1$ is centered, or 
\[
\varphi'(a_1ja_2)=\varphi'(j)\varphi(a_1a_2)=\varphi(a_1a_2),
\]
using the Leibniz-type formula of Definition 2.1. The last case gives (2). 
It remains to prove that if $j$ occurs more than once, then all terms vanish on the RHS of 
the Leibniz-type formula of Definition 2.1. Namely, there are two types of terms: 
those containing $\varphi'(a_k)$ for some $k$ and those containing 
$\varphi'(j)$. In the first case, there will always survive at least one $j$ under $\varphi$ since $j^p=j$ for any natural $p$ 
(since $j$ is an idempotent). Therefore, since $\varphi(j^p)=\varphi(j)=0$, such terms vanish by freeness condition (1) of Definition 2.1.
In the second case, after pulling out $\varphi'(j)$, we still have at least one $j$ under $\varphi$ and therefore 
this moment under $\varphi$ vanishes by the same freeness condition. Observe that this remains true after reduction 
$a_{k-1}a_{k+1}=a+c1$, where $a$ is centered and $c\in {\mathbb C}$. This proves (1) and thus the proof is complete.
\hfill $\blacksquare$\\

In the framework of Mingo and Tseng \cite{MT}, the functional $\varphi'$ with $\varphi'(1)=0$ leads to a normalized linear functional 
$\psi$ on ${\mathcal A}$ which has a nice relation to Boolean independence.

\begin{Corollary}
Under the assumptions of Proposition 2.2, let $\psi: {\mathcal A}\rightarrow {\mathbb C}$ be the linear functional defined by 
the formula 
\begin{equation}\label{eq:2.4}
\psi(a)=\varphi'(aj)
\end{equation}
for any $a\in {\mathcal A}$. Then it holds that
\begin{equation}\label{eq:2.5}
\psi(j^{\alpha}a_1j\cdots ja_{n}j^{\beta})=\psi(a_1)\cdots \psi(a_n)
\end{equation}
for any $a_1, \ldots, a_n\in {\mathcal B}$ and any $\alpha, \beta\in \{0,1\}$. Moreover, $\psi(1)=\psi(j)=1$ and 
$\psi(a)=\varphi(a)$ for any $a\in {\mathcal B}$. 
\end{Corollary}
\vspace{3pt}

\begin{Remark}
{\rm Certain important mixed moments w.r.t. $\psi$ can be expressed in terms of Boolean cumulants. Namely, 
it was shown in \cite{MT}, in Proposition 5.3, that 
\begin{equation}\label{eq:2.6}
\psi(a_1j^{\perp}a_{2}j^{\perp}\cdots j^{\perp}a_{n})=\beta_n(a_1, \ldots, a_n)
\end{equation}
where $a_1,\ldots, a_n\in {\mathcal B}$ and $j^{\perp}=1-j$, with $\beta_n(a_1, \ldots, a_n)$ denoting
the Boolean cumulants associated with moments under $\varphi$. An equivalent result was proved by the author in \cite{L3},
but we used partitioned `inverse Boolean cumulants'.
}
\end{Remark}
\vspace{3pt}
By Corollary 2.1, $\psi$ is an extension of $\varphi|_{\mathcal B}$ and 
the (non-unital) algebra generated by $j$ is Boolean independent against ${\mathcal B}$ with respect 
to $\psi$. These facts played a fundamental role in our approach to the unification of noncommutative independence 
(or, reduction of freeness and Boolean independence to tensor independence) in \cite{L1, L2}. We recall our definition below.

\begin{Definition}
{\rm Let $({\mathcal A}, \varphi)$ be a noncommutative probability space and let $p$ be an idempotent. 
By a {\it Boolean extension} of $\varphi$ we understand the functional $\widetilde{\varphi}$ on
\begin{equation}\label{eq:2.7}
\widetilde{\mathcal A}={\mathcal A}\star {\rm alg}(\{p\})
\end{equation}
given by the linear extension of 
\begin{equation}\label{eq:2.8}
\widetilde{\varphi}(p^{\alpha}a_1p\cdots pa_np^{\beta})
=
\varphi(a_1)\cdots \varphi(a_n)
\end{equation}
where $\alpha, \beta\in \{0,1\}$ and $a_1,\ldots , a_n\in {\mathcal A}$.
The noncommutative probability space $(\widetilde{\mathcal A}, \widetilde{\varphi})$ 
is called a {\it Boolean extension} of $({\mathcal A}, \varphi)$.
}
\end{Definition}

Therefore, if we regard $\psi$ as a functional defined on the free product 
${\mathcal B}\star {\rm alg}(j)$, the functional $\psi$ 
given by \cref{eq:2.4,eq:2.5} can be identified with our 
Boolean extension $\widetilde{\varphi}$, with the idempotent $j$ replaced by
$p$. Although the Boolean extension and the infinitesimal idempotent yield 
the same mixed moments, they arise from conceptually different constructions. 
In the Boolean extension framework, one algebraically adjoins an idempotent $p$ and defines an extended functional 
on the enlarged algebra. In contrast, in the setting of Mingo and Tseng \cite{MT} and
Shlyakhtenko \cite{Sh}, the idempotent $j$ is assumed to belong to an algebra and it is 
infinitesimally free with respect to a derivative functional $\varphi'$ and the functional $\psi$
is obtained from $\varphi'$ by a conditioning procedure given by \eqref{eq:2.4}. 
This correspondence allows us to reinterpret our construction in the infinitesimal framework. Since freeness
in our approach is obtained from a countable tensor product of Boolean extensions, the construction 
of freeness in \cite{L1} may thus be viewed as being consistent with a countable family of infinitesimally 
free objects while obtained through a different algebraic mechanism.

Apart from this conceptual connection, we use Boolean extensions of a continuously deformed family of functionals to construct a deformed tensor product functional. This framework allows us to apply derivations on tensor products and thereby define infinitesimal tensor product functionals associated with other notions of independence, whenever tensor product constructions of the same type are available. In particular, the canonical derivative of a tensor product functional is given by the sum over tensor sites of the derivative at a single site tensored with the undeformed functionals at all remaining sites.\\

\section{Motzkin paths and level return partitions}

In order to investigate the moments of orthogonal replicas, we need to
study Motzkin paths and their combinatorics. Let us begin with a standard definition 
of a Motzkin path. 

\begin{Definition}
{\rm By a {\it Motzkin path} of length $n-1$ we shall understand a lattice path 
\[
(v_1, \ldots v_n),
\]
where $v_k=(x_k,y_k)$ and $y_k\geq 0$ for each $k$, with the initial point $(0,0)$ and terminal point $(n-1,0)$, 
such that $v_{k}-v_{k-1}\in \{U, H, D\}$, where $U=(1,1),\,H=(1,0),\,D=(1,-1)$ are up, horizontal and down steps, respectively, 
for each $k=2,\ldots, n$. We identify this path with the corresponding step-word built from letters $U,H,D$.}
\end{Definition}

In our approach, we will use still another description of 
Motzkin paths, which encodes the $y$-coordinates of the vertices.
These coordinates, as defined above, are nonnegative 
integers, but in our approach to Motzkin 
paths we find it convenient to use only positive integers. 
Therefore, we will identify Motzkin paths not only with step-words, but also
with level-words built from integers $j_k=y_k+1$, where $k=1, \ldots, n$. Since 
in our definition of cumulants in \cite{L5} we also used subpaths of Motzkin paths 
(starting and ending with any $j_1=j_n\in {\mathbb N}$), we distinguished
Motzkin words starting and ending with $j_1=j_n=1$ and called these words {\it reduced}.

\begin{Definition}
{\rm By a {\it reduced Motzkin word} of length $n\in \mathbb{N}$ we 
understand a word in letters from the alphabet 
${\mathbb N}$ of the form
\[
w=j_1\cdots j_n
\]
where $j_1, \ldots, j_n\in \mathbb{N}$, $j_1=j_n=1$ and $j_{i}-j_{i-1}\in \{-1,0,1\}$ for any $i=2, \ldots, n$.
The set of reduced Motzkin words of length $n$ will be denoted by $\mathpzc{M}_{n}$. Then
\[
\mathpzc{M}=\bigcup_{n=1}^{\infty}\mathpzc{M}_{n}
\] 
will stand for the set of all reduced Motzkin words.}
\end{Definition}

Note that a reduced Motzkin word of length $n$ is associated with a Motzkin path of length $n-1$, while the word $w=1$ is associated with the empty path. Occasionally, we will use the {\it height} of $w$, defined by 
\[
h(w)={\rm max}\{j_1,\ldots, j_n\}.
\]
This differs from the geometric height of the associated path since that one
is equal to $h(w)-1$. To avoid duplicating terminology and notation for paths and words, 
we will work exclusively with Motzkin words $w$ in what follows. 
When referring to objects (such as partitions or functionals) associated with $w$, 
we will nevertheless say that they are associated with the Motzkin path $w$.
This mild abuse of language should cause no confusion and is convenient, since Motzkin paths have a natural geometric interpretation and admit analytic features (such as local maxima), whereas restricting attention to words alone would be unnecessarily limiting.

Let us introduce the key new combinatorial notion of a {\it level return partition}
associated with $w\in \mathpzc{M}_{n}$. The definition refers 
to the heights ({\it colors} or {\it levels}) of the vertices, 
namely $j_1, \ldots, j_n$. First we define blocks of this partition. 
It it convenient to require that two consecutive numbers do not
belong to the same block (we call it the `gap condition'). 
This anticipates the fact that the neighboring variables in a 
moment associated with a Motzkin path will be assumed to have different labels, 
although at this point we would like to define the blocks 
of a level return partition knowing only $w=j_1\cdots j_n$.

\begin{Definition}
{\rm For each $w=j_1\cdots j_n\in \mathpzc{M}_{n}$ and each 
$j\in \{j_1, \ldots, j_n\}$, let 
\[
S_{j}=\{k\in [n]: j_k=j\}=\{k_1<\cdots < k_m\},\;\;{\rm where}\;m\geq 1.
\]
Any subset $V=\{k_{p}<\cdots <k_{q}\}\subseteq S_j$, where $p\leq q$, 
will be called a {\it level return block of color} $j$ if it
satisfies the following conditions:
\begin{enumerate}
\item 
{\it gap condition}: 
\[
k_{s+1}-k_{s}>1,
\]
\item
{\it return condition}:
\[
k_{s}<r<k_{s+1}\Longrightarrow  j_r> j,
\]
\end{enumerate}
for all $p\leq s <q$.
In particular, if $p=q$, then $V$ is a singleton block of color $j$.
Elements of a level return block are called {\it level return points} 
and the corresponding vertices of the path are called {\it level return vertices}.}
\end{Definition}

\begin{Definition}
{\rm A level return block $V\subseteq S_j$ will be called {\it maximal} if it is not 
properly contained in any other level return block of color $j$. 
The partition consisting of maximal level return blocks will be called 
a {\it level return partition associated with $w$} and will be denoted $\pi(w)$.}
\end{Definition}

\begin{Definition}
{\rm 
Let $w=j_1\cdots j_n\in \mathpzc{M}_{n}$.
By an {\it excursion above level $j$} (or, simply, an {\it excursion}) in $w$ 
we understand a subword $w'=j_{s+1}\cdots j_{t-1}$, such that 
$w=j_1\cdots j_{s}w'j_{t}\cdots j_n$, where $j_s=j_t=j$ and $j_k>j$ for all points $s<k<t$.
Geometrically, an excursion above level $j$ is a subpath between two 
level return vertices lying at level $j$ that stays above level $j$.
From now on we will identify `color' with `level' for vertices, blocks and excursions. 
}
\end{Definition}

\begin{Lemma}
For every \(w\in\mathpzc{M}_n\) there exists a unique level return partition $\pi(w)$
and this partition is noncrossing.
\end{Lemma}
{\it Proof.}
For any $j$, let $S_j=\{k_1<\cdots <k_m\}$.
Every element $t\in S_j$ belongs to a maximal subsequence of the form
\[
k_{p}< \cdots < k_{q}
\]
where $p\leq q$, such that each pair $(k_{s},k_{s+1})$ consists of non-consecutive indices 
which satisfy the return condition
for each $p\leq s< q$ (if $p=q$ then $\{t\}$ is a singleton).
Every element of $S_{j}$ belongs to a unique maximal subsequence of the above form which is a 
maximal level return block.
Repeating this for all $j$ produces a unique family of disjoint maximal
level return blocks. This proves uniqueness. 
To show that $\pi(w)$ is noncrossing, suppose two distinct maximal blocks \(V=\{s_1<\dots<s_r\}\) 
at level \(j\) and \(U=\{t_1<\dots<t_s\}\) at level \(j'\) cross, so for some indices
\[
s_i<t_k<s_{i+1}<t_{k+1}.
\]
By the return condition, every index between \(s_i\) and \(s_{i+1}\) has color \(\ge j\), hence 
\(j_{t_k}\ge j\); similarly every index between \(t_k\) and \(t_{k+1}\) has color \(\ge j'\), hence \(j_{s_{i+1}}\ge j'\). Thus \(j\ge j'\) and \(j'\ge j\), so \(j=j'\). But then the union \(V\cup U\) satisfies the return condition,
contradicting maximality of \(V\) and \(U\). Therefore no crossing can occur, and \(\pi(w)\) is noncrossing.
\hfill $\blacksquare$\\

\begin{Definition}
{\rm The level return partition $\pi(w)$ is {\it adapted} to the sequence of labels $(i_1, \ldots, i_n)$, where $i_1\neq \ldots \neq i_n$, 
if
\begin{enumerate}
\item for any two level return points $p$ and $q$ belonging to the same block $V\in \pi$
it holds that $i_p=i_q$ ({\it label uniformity}),
\vspace{3pt}
\item
for any two consecutive level return points $p$ and $q$ belonging to the same 
block $V\in \pi$ at level $j$ and any block $V'\in \pi$ at level $j+1$ contained
in the excursion between $j_p$ and $j_q$, the following condition holds ({\it nested alternation}):
\[
i_p\neq i_{p'}\;\;{\rm for\;any}\;p\in V\;\;{\rm and\;any}\;p'\in V'.
\]
The labels of adjacent blocks of such a partition must therefore alternate.
\end{enumerate}
}
\end{Definition}
\vspace{3pt}

\begin{Example}
{\rm 
Anticipating that our combinatorics will be applied to infinitesimal freeness, observe that 
in the formula for centered moments of free random variables we have $i_1\neq \ldots \neq i_n$
and 
\[
i_1=i_n,\;\; i_2=i_{n-1},\;\;\ldots, \;\;i_{m-1}=i_{m+1}
\]
for odd $n=2m-1$. If we imagine that this sequence may correspond to a {\it pyramid word} (see Fig. 2) 
\[
w=12\cdots (m-1)m(m+1)\cdots n
\]
then points $1,n$ are level return points and the associated partition 
\[
\pi(w)=\{\{1,n\}, \{2, n-1\}, \ldots, \{m-1,m+1\}, \{m\}\}
\] 
is the level return partition adapted to the tuple of labels $\ell=(i_1, \ldots, i_n)$.
}
\end{Example}

\begin{Example}
{\rm Consider three examples of Motzkin paths with various localizations of horizontal steps 
and the corresponding level return partitions shown in Fig. 1.
\begin{enumerate}
\item
For the word $w_1=123332112121$, the height is $3$, thus the returns occur at $3$ levels: $1,2,3$. Namely,
\[
S_{1}=\{1,7,8,10,12\}, \;\;S_{2}=\{2,6,9,11\}, \;\;S_{3}=\{3,4,5\}.
\]
The associated level return partition is 
\[
\pi(w_1)=\{\{1,7\}, \{2,6\},\{3\},\{4\},\{5\},\{8,10,12\}, \{9\},\{11\} \}.
\]
Here the refinement of $\{S_1,S_2,S_3\}$ is determined by the gap condition and by whether the path
goes below a given level between successive returns.
Adaptedness to $(i_1,\ldots,i_{12})$ means that labels are constant on level return blocks
($i_1=i_7$, $i_2=i_6$, $i_8=i_{10}=i_{12}$), and the nested alternation condition yields
$i_2\neq i_4$, in addition to $i_1\neq\cdots\neq i_{12}$.
\item
For $w_2=112323223211$, we have
\[
S_{1}=\{1,2,11,12\}, \;
S_{2}=\{3,5,7,8,10\}, \;
S_{3}=\{4\}, \;
S_{4}=\{6\}, \;
S_{5}=\{9\},
\]
and
\[
\pi(w_2)=\{\{1\}, \{2,11\},\{3,5,7\},\{4\},\{6\},\{8,10\}, \{9\},\{12\} \}.
\]
Adaptedness requires constant labels on blocks
($i_2=i_{11}$, $i_3=i_5=i_7$, $i_8=i_{10}$), and the nested alternation condition imposes
no further restrictions beyond $i_1\neq\cdots\neq i_{12}$.
\item
For $w_3=123432334321$, we have
\[
S_{1}=\{1,12\}, \;
S_{2}=\{2,6,11\}, \;
S_{3}=\{3,5,7,8,10\}, \;
S_{4}=\{4,9\},
\]
and the level return partition is
\[
\pi(w_3)=\{\{1,12\}, \{2,6,11\},\{3,5\},\{4\},\{7\},\{8,10\},\{9\}\}.
\]
Adaptedness means that labels are constant on level return blocks
($i_1=i_{12}$, $i_2=i_6=i_{11}$, $i_8=i_{10}$), with no additional constraints
besides $i_1\neq\cdots\neq i_{12}$.
\end{enumerate}}
\end{Example}

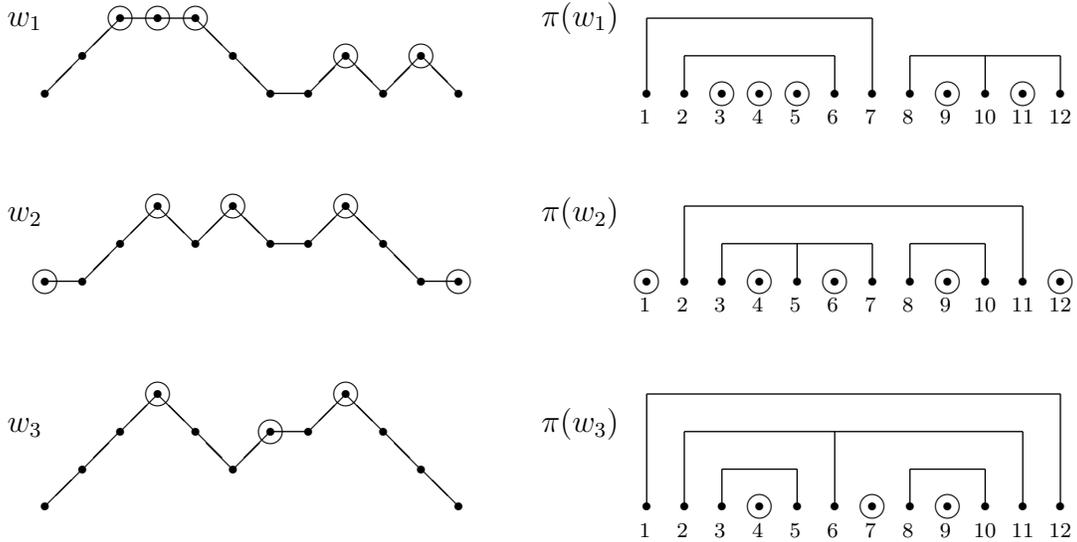
\begin{figure}
\unitlength=1mm
\special{em:linewidth 1pt}
\linethickness{0.5pt}
\begin{picture}(180.00,85.00)(-22.00,-10.00)

%%%%%%%%%%%%%%%path %%%%%%%%%%%%%%%
\put(-10.00,64.00){$w_1$}
\put(-10.00,38.00){$w_2$}
\put(-10.00,10.00){$w_3$}

\put(61.00,64.00){$\pi(w_1)$}
\put(61.00,38.00){$\pi(w_2)$}
\put(61.00,10.00){$\pi(w_3)$}

%%%%%%%%%%%%%%%%% path 1%%%%%%%%%%%%%%%%%%%%%%%%
\put(-5.00,55.00){\line(1,1){5.00}}
\put(0.00,60.00){\line(1,1){5.00}}
\put(5.00,65.00){\line(1,0){5.00}}
\put(10.00,65.00){\line(1,0){5.00}}
\put(15.00,65.00){\line(1,-1){5.00}}
\put(20.00,60.00){\line(1,-1){5.00}}
\put(25.00,55.00){\line(1,0){5.00}}
\put(30.00,55.00){\line(1,1){5.00}}
\put(35.00,60.00){\line(1,-1){5.00}}
\put(40.00,55.00){\line(1,1){5.00}}
\put(45.00,60.00){\line(1,-1){5.00}}

\put(-5.00,55.00){\circle*{1.00}}
\put(0.00,60.00){\circle*{1.00}}
\put(5.00,65.00){\circle*{1.00}}
\put(5.00,65.00){\circle{3.00}}
\put(10.00,65.00){\circle*{1.00}}
\put(10.00,65.00){\circle{3.00}}
\put(15.00,65.00){\circle*{1.00}}
\put(15.00,65.00){\circle{3.00}}
\put(20.00,60.00){\circle*{1.00}}
\put(25.00,55.00){\circle*{1.00}}
\put(30.00,55.00){\circle*{1.00}}
\put(35.00,60.00){\circle*{1.00}}
\put(35.00,60.00){\circle{3.00}}
\put(40.00,55.00){\circle*{1.00}}
\put(45.00,60.00){\circle*{1.00}}
\put(45.00,60.00){\circle{3.00}}
\put(50.00,55.00){\circle*{1.00}}

%%%%partition 1 %%%%%%%%%%%%%%%%%%%%%%

\put(75.00,55.00){\circle*{1.00}}
\put(80.00,55.00){\circle*{1.00}}
\put(85.00,55.00){\circle*{1.00}}
\put(85.00,55.00){\circle{3.00}}
\put(90.00,55.00){\circle*{1.00}}
\put(90.00,55.00){\circle{3.00}}
\put(95.00,55.00){\circle*{1.00}}
\put(95.00,55.00){\circle{3.00}}
\put(100.00,55.00){\circle*{1.00}}
\put(105.00,55.00){\circle*{1.00}}
\put(110.00,55.00){\circle*{1.00}}
\put(115.00,55.00){\circle*{1.00}}
\put(115.00,55.00){\circle{3.00}}
\put(120.00,55.00){\circle*{1.00}}
\put(125.00,55.00){\circle*{1.00}}
\put(125.00,55.00){\circle{3.00}}
\put(130.00,55.00){\circle*{1.00}}

\put(74.00,51.00){$\scriptstyle{1}$}
\put(79.00,51.00){$\scriptstyle{2}$}
\put(84.00,51.00){$\scriptstyle{3}$}
\put(89.00,51.00){$\scriptstyle{4}$}
\put(94.00,51.00){$\scriptstyle{5}$}
\put(99.00,51.00){$\scriptstyle{6}$}
\put(104.00,51.00){$\scriptstyle{7}$}
\put(109.00,51.00){$\scriptstyle{8}$}
\put(114.00,51.00){$\scriptstyle{9}$}
\put(118.50,51.00){$\scriptstyle{10}$}
\put(123.50,51.00){$\scriptstyle{11}$}
\put(128.50,51.00){$\scriptstyle{12}$}

\put(75.00,55.00){\line(0,1){10}}
\put(80.00,55.00){\line(0,1){5}}
\put(100.00,55.00){\line(0,1){5}}
\put(105.00,55.00){\line(0,1){10}}
\put(110.00,55.00){\line(0,1){5}}
\put(120.00,55.00){\line(0,1){5}}
\put(130.00,55.00){\line(0,1){5}}
\put(75.00,65.00){\line(1,0){30}}
\put(80.00,60.00){\line(1,0){20}}
\put(110.00,60.00){\line(1,0){20}}
%\put(125.00,60.00){\line(1,0){5}}

%%%%%%%%%%%%%% path 2 %%%%%%%%%%%%%%%%%%%%%%%%%%%%%%
\put(-5.00,30.00){\line(1,0){5.00}}
\put(0.00,30.00){\line(1,1){5.00}}
\put(5.00,35.00){\line(1,1){5.00}}
\put(10.00,40.00){\line(1,-1){5.00}}
\put(15.00,35.00){\line(1,1){5.00}}
\put(20.00,40.00){\line(1,-1){5.00}}
\put(25.00,35.00){\line(1,0){5.00}}
\put(30.00,35.00){\line(1,1){5.00}}
\put(35.00,40.00){\line(1,-1){5.00}}
\put(40.00,35.00){\line(1,-1){5.00}}
\put(45.00,30.00){\line(1,0){5.00}}

\put(-5.00,30.00){\circle*{1.00}}
\put(-5.00,30.00){\circle{3.00}}
\put(0.00,30.00){\circle*{1.00}}
\put(5.00,35.00){\circle*{1.00}}
\put(10.00,40.00){\circle*{1.00}}
\put(10.00,40.00){\circle{3.00}}
\put(15.00,35.00){\circle*{1.00}}
\put(20.00,40.00){\circle*{1.00}}
\put(20.00,40.00){\circle{3.00}}
\put(25.00,35.00){\circle*{1.00}}
\put(30.00,35.00){\circle*{1.00}}
\put(35.00,40.00){\circle*{1.00}}
\put(35.00,40.00){\circle{3.00}}
\put(40.00,35.00){\circle*{1.00}}
\put(45.00,30.00){\circle*{1.00}}
\put(50.00,30.00){\circle*{1.00}}
\put(50.00,30.00){\circle{3.00}}

%%%%partition 2 %%%%%%%%%%%%%%%%%%%%%%

\put(75.00,30.00){\circle*{1.00}}
\put(75.00,30.00){\circle{3.00}}
\put(80.00,30.00){\circle*{1.00}}
\put(85.00,30.00){\circle*{1.00}}
\put(90.00,30.00){\circle*{1.00}}
\put(90.00,30.00){\circle{3.00}}
\put(95.00,30.00){\circle*{1.00}}
\put(100.00,30.00){\circle*{1.00}}
\put(100.00,30.00){\circle{3.00}}
\put(105.00,30.00){\circle*{1.00}}
\put(110.00,30.00){\circle*{1.00}}
\put(115.00,30.00){\circle*{1.00}}
\put(115.00,30.00){\circle{3.00}}
\put(120.00,30.00){\circle*{1.00}}
\put(125.00,30.00){\circle*{1.00}}
\put(130.00,30.00){\circle*{1.00}}
\put(130.00,30.00){\circle{3.00}}

\put(80.00,30.00){\line(0,1){10}}
\put(85.00,30.00){\line(0,1){5}}
\put(95.00,30.00){\line(0,1){5}}
\put(105.00,30.00){\line(0,1){5}}
\put(110.00,30.00){\line(0,1){5}}
\put(120.00,30.00){\line(0,1){5}}
\put(125.00,30.00){\line(0,1){10}}
\put(80.00,40.00){\line(1,0){45}}
\put(85.00,35.00){\line(1,0){20}}
\put(110.00,35.00){\line(1,0){10}}

\put(74.00,26.00){$\scriptstyle{1}$}
\put(79.00,26.00){$\scriptstyle{2}$}
\put(84.00,26.00){$\scriptstyle{3}$}
\put(89.00,26.00){$\scriptstyle{4}$}
\put(94.00,26.00){$\scriptstyle{5}$}
\put(99.00,26.00){$\scriptstyle{6}$}
\put(104.00,26.00){$\scriptstyle{7}$}
\put(109.00,26.00){$\scriptstyle{8}$}
\put(114.00,26.00){$\scriptstyle{9}$}
\put(118.50,26.00){$\scriptstyle{10}$}
\put(123.50,26.00){$\scriptstyle{11}$}
\put(128.50,26.00){$\scriptstyle{12}$}

%%%%%%%%%%%%%%% path 3 %%%%%%%%%%%%%%%%%%%%%%%%%%%%%

\put(-5.00,0.00){\line(1,1){5.00}}
\put(0.00,5.00){\line(1,1){5.00}}
\put(5.00,10.00){\line(1,1){5.00}}
\put(10.00,15.00){\line(1,-1){5.00}}
\put(15.00,10.00){\line(1,-1){5.00}}
\put(20.00,5.00){\line(1,1){5.00}}
\put(25.00,10.00){\line(1,0){5.00}}
\put(30.00,10.00){\line(1,1){5.00}}
\put(35.00,15.00){\line(1,-1){5.00}}
\put(40.00,10.00){\line(1,-1){5.00}}
\put(45.00,5.00){\line(1,-1){5.00}}

\put(-5.00,0.00){\circle*{1.00}}
\put(0.00,5.00){\circle*{1.00}}
\put(5.00,10.00){\circle*{1.00}}
\put(10.00,15.00){\circle*{1.00}}
\put(10.00,15.00){\circle{3.00}}
\put(15.00,10.00){\circle*{1.00}}
\put(20.00,5.00){\circle*{1.00}}
\put(25.00,10.00){\circle*{1.00}}
\put(25.00,10.00){\circle{3.00}}
\put(30.00,10.00){\circle*{1.00}}
\put(35.00,15.00){\circle*{1.00}}
\put(35.00,15.00){\circle{3.00}}
\put(40.00,10.00){\circle*{1.00}}
\put(45.00,5.00){\circle*{1.00}}
\put(50.00,0.00){\circle*{1.00}}

%%%% partition 3 %%%%%%%%%%%%%%%%%%%%%%

\put(75.00,0.00){\circle*{1.00}}
\put(80.00,0.00){\circle*{1.00}}
\put(85.00,0.00){\circle*{1.00}}
\put(90.00,0.00){\circle*{1.00}}
\put(90.00,0.00){\circle{3.00}}
\put(95.00,0.00){\circle*{1.00}}
\put(100.00,0.00){\circle*{1.00}}
\put(105.00,0.00){\circle*{1.00}}
\put(105.00,0.00){\circle{3.00}}
\put(110.00,0.00){\circle*{1.00}}
\put(115.00,0.00){\circle*{1.00}}
\put(115.00,0.00){\circle{3.00}}
\put(120.00,0.00){\circle*{1.00}}
\put(125.00,0.00){\circle*{1.00}}
\put(130.00,0.00){\circle*{1.00}}

\put(74.00,-4.00){$\scriptstyle{1}$}
\put(79.00,-4.00){$\scriptstyle{2}$}
\put(84.00,-4.00){$\scriptstyle{3}$}
\put(89.00,-4.00){$\scriptstyle{4}$}
\put(94.00,-4.00){$\scriptstyle{5}$}
\put(99.00,-4.00){$\scriptstyle{6}$}
\put(104.00,-4.00){$\scriptstyle{7}$}
\put(109.00,-4.00){$\scriptstyle{8}$}
\put(114.00,-4.00){$\scriptstyle{9}$}
\put(118.50,-4.00){$\scriptstyle{10}$}
\put(123.50,-4.00){$\scriptstyle{11}$}
\put(128.50,-4.00){$\scriptstyle{12}$}

\put(75.00,0.00){\line(0,1){15}}
\put(80.00,0.00){\line(0,1){10}}
\put(85.00,0.00){\line(0,1){5}}
\put(95.00,0.00){\line(0,1){5}}
\put(100.00,0.00){\line(0,1){10}}
\put(110.00,0.00){\line(0,1){5}}
\put(120.00,0.00){\line(0,1){5}}
\put(125.00,0.00){\line(0,1){10}}
\put(130.00,0.00){\line(0,1){15}}

\put(75.00,15.00){\line(1,0){55}}
\put(80.00,10.00){\line(1,0){45}}
\put(85.00,5.00){\line(1,0){10}}
\put(110.00,5.00){\line(1,0){10}}

\end{picture}
\caption{Three Motzkin paths and the associated level return partitions. 
Local maxima of paths and the corresponding singletons in level return partitions 
are marked with hollow circles.}
\end{figure}

An important property of level return partitions is that singleton blocks of $\pi(w)$ 
correspond to {\it local maxima} of the corrresponding path.

\begin{Lemma}
The set $V=\{k\}$, where $k\in [n]$, is a singleton block of $\pi(w)$ if and only if
the corresponding Motzkin path has a local maximum at $(k,j_k)$, that is 
\[
j_{k-1}\leq j_k\geq j_{k+1}
\] 
for the corresponding word $w=j_1\cdots j_n$, with the understanding that if 
$k=1$ or $k=n$, then only  the right or left inequality holds, respectively.
\end{Lemma}
{\it Proof.}
The set $\{k\}$ is a singleton block of color $j_k$ of $\pi(w)$ if and only if
$k$ is not connected to the largest index $k'<k$ with $j_{k'}=j_k$
nor to the smallest index $k''>k$ with $j_{k''}=j_k$. 
By the definition of the construction, this occurs if and only if, on both sides of $k$, 
at least one of the gap or return condition fails. 
Translating these failures into inequalities between consecutive heights, we obtain
\[
j_{k-1}\leq j_k\;\;\; {\rm and}\;\;\;j_{k+1}\leq j_k
\]
with the obvious modifications when $k=1$ or $k=n$. Equivalently, the 
Motzkin path satisfies $j_{k-1}\leq j_k \geq j_{k+1}$ , so the path has a local maximum at 
$(k,j_k)$. This completes the proof. 
\hfill $\blacksquare$\\

\section{Infinitesimal deformations of Motzkin functionals}

We showed in \cite{L4} that the mixed moments of free random variables 
decompose in terms of the mixed moments of variables called orthogonal replicas
with respect to a tensor product functional, as in \cref{eq:1.3,eq:1.4}.

In the case of two algebras, these variables have a simple form given below \cite{L5} and 
we would like to write it explicitly before treating the case of an arbitrary index set $I$.
The definition is based on the concept of a Boolean extension of a noncommutative probability space
(see Definition 2.2). 

\begin{Definition}
{\rm 
Let $({\mathcal A}_{1},\varphi_{1})$ and $({\mathcal A}_{2}, \varphi_{2})$ be noncommutative probability spaces and
let $(\widetilde{\mathcal A}_{1}, \widetilde{\varphi}_{1})$ and $(\widetilde{\mathcal A}_{2}, \widetilde{\varphi}_{2})$ 
be their Boolean extensions (with idempotents $p_1$ and $p_2$, respectively). Define infinite tensor products
{
\setlength{\belowdisplayshortskip}{\baselineskip}
\begin{equation}\label{eq:4.1}
{\mathcal A}_{\otimes}:=
\widetilde{\mathcal A}_{1}^{\otimes \infty}\otimes \widetilde{\mathcal A}_{2}^{\otimes\infty}\;\;\;{\rm and}\;\;\;
\Phi_{\otimes}:=\widetilde{\varphi}_{1}^{\otimes\infty}\otimes \widetilde{\varphi}_{2}^{\otimes\infty}.
\end{equation}
}%
The elements of ${\mathcal A}_{\otimes}$ of the form\\
\begin{align}\label{eq:4.2}
a(j)&=(\cdots \otimes 1_1 \otimes \underbrace{a}_{{\rm site}\;j} \otimes 1_1\otimes \cdots )\otimes 
(\cdots \otimes 1_2\otimes \underbrace{p_2^{\perp}}_{{\rm site}\;j-1} \otimes\; p_2 \otimes \cdots),\\
\label{eq:4.3}
a'(j)&=(\cdots \otimes 1_1 \otimes \underbrace{p_1^{\perp}}_{{\rm site}\;j-1} \otimes p_1\otimes \cdots )\otimes 
(\cdots \otimes 1_2\otimes \underbrace{a'}_{{\rm site}\;j} \otimes\; 1_2 \otimes \cdots)
\end{align}
where $a\in \mathcal{A}_{1}$, $a'\in \mathcal{A}_2$, and 
$1_k$ is the internal unit in ${\mathcal A}_{k}$ for $k=1,2$,
will be called {\it orthogonal replicas} of $a$ and $a'$ of {\it color} $j\in \mathbb{N}$ and 
{\it labels} $1$ and $2$, respectively. We understand that if $j=1$, 
then the expressions marked with `site $j-1$' do not appear.
}
\end{Definition}

For an arbitrary index set $I$, the orthogonal replicas of color $j$ 
are not of the simple form given by \cref{eq:4.2,eq:4.3}. Instead of local $p_1^{\perp}$ and $p_2^{\perp}$ 
at sites of color $j-1$ and labels different from the given $i$ 
we have global orthogonal projections. Let us now give the general definition.
Thus, let
{
\setlength{\belowdisplayshortskip}{\baselineskip}
\begin{equation}\label{eq:4.4}
{\mathcal A}_{\otimes}:=\bigotimes_{i\in I}
\widetilde{\mathcal A}_{i}^{\otimes \infty}\;\;\;{\rm and}\;\;\;
\Phi_{\otimes}:=\bigotimes_{i\in I}\widetilde{\varphi}_{i}^{\otimes\infty}
\end{equation}
}%
be the infinite tensor products for an arbitrary index set $I$.
\begin{Definition}
{\rm 
By {\it orthogonal replicas} of variables $a_{k}\in {\mathcal A}_{i_{k}}$ we understand
elements of the tensor product algebra ${\mathcal A}_{\otimes}$ of the form 
\begin{equation}\label{eq:4.5}
a_{k}(j_k)=(a_{k})_{i_k,j_k}\otimes (P_{i_k,j_k}-P_{i_k,j_k-1}),
\end{equation}
where 
\begin{equation}\label{eq:4.6}
P_{i_k,j_k}=\bigotimes_{i\neq i_k, j\geq j_k}(p_{i})_{i,j}.
\end{equation}
Here, $(a_{k})_{i_k,j_k}$ indicates that $a_k$ is at site of label $i_k$ and color $j_k$ and 
at the remaining tensor positions of label $i_k$ we have an internal unit $1_{i_k}$. 
Similarly, in $P_{i_k,j_k}$ the idempotent $p_{i}$ is at positions of label $i\neq i_k$ 
and color $j\geq j_k$ and at the remaining positions
of color $j$ and label $i\neq i_k$ (omitted in the equation) there is an internal unit $1_i$.}
\end{Definition}

The noncommutative probability space $({\mathcal A}, \Phi)$, where ${\mathcal A}$ is the unital 
subalgebra of ${\mathcal A}_{\otimes}$ defined in \eqref{eq:4.4}
generated by orthogonal replicas and $\Phi=\Phi_{\otimes}|_{\mathcal A}$, 
is called the {\it orthogonal replica space}. In view of the infinitesimal structure 
inherited from \cite{MT}, we can view this space as a space equipped with a countable number of 
infinitesimal idempotents. For our purposes, $({\mathcal A}, \Phi)$ is sufficiently large, but 
we can also use the much larger $({\mathcal A}_{\otimes}, \Phi_{\otimes})$ for moment computations. 
The functional $\Phi$ is a tool to define a family of 
multilinear Motzkin functionals introduced in \cite{L4}, where we denoted 
them by $\psi(w)$. We will use a different notation in this paper, namely $\Phi_{w}$, since in 
what follows $w$ is typically fixed and emphasis is put on the arguments
$a_1, \ldots, a_n$. This notation avoids nested parentheses and highlights the role 
of $w$ as an index.

\begin{Definition}
{\rm 
By {\it Motzkin functionals} we understand the multilinear moment functionals $\Phi_w$, where $w=j_1\cdots j_n\in \mathpzc{M}_{n}$, 
defined as multilinear extensions of
{
\setlength{\belowdisplayshortskip}{\baselineskip}
\begin{equation}\label{eq:4.7}
\Phi_{w}(a_1, \ldots, a_n)=\Phi(a_1(j_1)\cdots a_n(j_n))
\end{equation}
}%
where $a_1(j_1), \ldots, a_n(j_n)$ are the orthogonal replicas of $a_1\in \mathcal{A}_{i_1}, \ldots, a_n\in \mathcal{A}_{i_n}$.}
\end{Definition}

Let us state the decomposition theorem of \cite{L4} in terms of Motzkin functionals. 
This decomposition expresses the free product of functionals 
in terms of Motzkin functionals. For a family of unital algebras $({\mathcal A}_{i})_{i\in I}$ 
denote by $\star_{i\in I}{\mathcal A}_i$ and $*_{i\in I}{\mathcal A}_{i}$ 
their free products with and without identification of units, respectively. 

\begin{Theorem}\cite{L4}
Let $({\mathcal A}_{i},\varphi_i)_{i\in I}$ be noncommutative probability spaces. 
It holds that
\begin{equation}\label{eq:4.8}
\varphi (a_1\cdots a_n)=
\sum_{w=j_1\cdots j_n\in \mathpzc{M}_{n}}
\Phi_{w}(a_{1}, \ldots, a_n)
\end{equation}
where $\varphi=(\star_{i\in I} \varphi_i)\circ \tau$
and $\tau: *_{i\in I}{\mathcal A}_i\rightarrow \star_{i\in I}{\mathcal A}_{i}$ is
the unit identification map, whenever $a_1\in \mathcal{A}_{i_1}, \ldots, a_n\in \mathcal{A}_{i_n}$ 
and $i_1\neq \cdots \neq i_n$, $n\in \mathbb{N}$.
\end{Theorem}

In order to study infinitesimal freeness (of any order) we
need to introduce a continuous deformation of the family $(\varphi_{i})$. 

\begin{Definition}
{\rm Let $(\varphi_{i})_{i\in I}$ be a family of normalized linear functionals 
on the family of unital algebras $({\mathcal A}_i)_{i\in I}$, respectively. Define a family of functionals 
of the form
\begin{equation}\label{eq:4.9}
\varphi_{i,t}=\varphi_{i}+t\varphi_{i}^{(1)}+ \frac{t^2}{2!}\varphi_{i}^{(2)} + \ldots +
\frac{t^{m}}{m!}\varphi_{i}^{(m)} + o(t^{m})
\end{equation}
where $i\in I$, $(\varphi_{i}^{(k)})_{k\geq 1}$ are additional linear functionals such that $\varphi_{i}^{(k)}(1_{i})=0$ for any $i\in I$ and any 
$1\leq k \leq m$, where $m\in {\mathbb N}$ and $t\in (-\epsilon, \epsilon)$. 
We will often write $\varphi'_{i}=\varphi_{i}^{(1)}$.}
\end{Definition}

The corresponding mixed moments and Boolean cumulants of order $n$ for 
elements $a_1\in {\mathcal A}_{i_1}, \ldots, a_n\in {\mathcal A}_{i_n}$ will be denoted
\[
\varphi_{i,t}(a_1\cdots a_n)\;\;\;{\rm and}\;\;\; 
\beta_{n,t}(a_1, \ldots, a_n),
\]
respectively. In our notation, we include the algebra label $i$ in the moments, 
but omit it in the case of Boolean cumulants in order to avoid notation with three indices, 
keeping the order $n$ of the cumulant as it is done by most authors. 
Therefore, the cumulant recognizes to which 
algebra its arguments belong.
\vspace{3pt}
\begin{Example}
{\rm The deformed Boolean cumulants of lowest orders are given by
\begin{align*}
\beta_{1,t}(a_1)&=\varphi_{i,t}(a_1),\\
\beta_{2,t}(a_1,a_2)&=\varphi_{i,t}(a_1a_2)-\varphi_{i,t}(a_1)\varphi_{i,t}(a_2),\\
\beta_{3,t}(a_1,a_2,a_3)&=\varphi_{i,t}(a_1a_2a_3)-\varphi_{i,t}(a_1a_2)
\varphi_{i,t}(a_3)-\varphi_{i,t}(a_1)\varphi_{i,t}(a_2a_3)\\
&\quad +\varphi_{i,t}(a_1)\varphi_{i,t}(a_2)\varphi_{i,t}(a_3),
\end{align*}
whenever $a_1,a_2,a_3\in {\mathcal A}_{i}$.
Their derivatives with respect to $t$ at $t=0$ are
\begin{align*}
\beta_{1}'(a_1)&=\varphi_{i}'(a_1),\\
\beta_{2}'(a_1,a_2)&=\varphi_{i}'(a_1a_2)-\varphi_{i}'(a_1)\varphi_{i}(a_2)
-\varphi_{i}(a_1)\varphi_{i}'(a_2),\\
\beta_{3}'(a_1,a_2,a_3)&=\varphi_{i}'(a_1a_2a_3)-\varphi_{i}'(a_1a_2)
\varphi_{i}(a_3)-\varphi_{i}(a_1a_2)\varphi_{i}'(a_3)\\
&\quad -\varphi_{i}'(a_1)\varphi_{i}(a_2a_3)-\varphi_{i}(a_1)\varphi_{i}'(a_2a_3)+
\varphi_{i}(a_1)'\varphi_{i}(a_2)\varphi_{i}(a_3)\\
&\quad +\varphi_{i}(a_1)\varphi_{i}'(a_2)\varphi_{i}(a_3)
+\varphi_{i}(a_1)\varphi_{i}(a_2)\varphi_{i}'(a_3).
\end{align*}
Therefore, also in the case of $\varphi_i'$ and $\beta_n'$ we use the same convention that 
$\varphi_i'$ shows the label of the functional, whereas $\beta_n$ shows the order of the cumulant.
}
\end{Example}
\vspace{3pt}

The construction relies on two local infinitesimal operations: a continuous deformation with parameter $t$, 
given by \eqref{eq:4.9}, and an operatorial one involving an infinitesimal idempotent $p_i$, 
given by \eqref{eq:2.8}:
\[
\varphi_{i}\xrightarrow {\text{deform}}
\varphi_{i,t}\xrightarrow {\text{extend}}
\widetilde{\varphi}_{i,t}
\]
for all $i\in I$ and $t$. Using these operations, we 
define a global continuous deformation of
$\Phi_{\otimes}$ and $\Phi$ by
\begin{equation}\label{eq:4.10}
\Phi_{\otimes, t}=\bigotimes_{i\in I}\widetilde{\varphi}_{i,t}^{\otimes \infty}\;\;\;{\rm and}\;\;\;\Phi_{t}=\Phi_{\otimes, t}|_{\mathcal A},
\end{equation}
respectively, thus $\Phi_{t}$ is the restriction of $\Phi_{\otimes}$ to 
the orthogonal replica space. We slightly abuse notation, since we also consider multilinear functionals $\Phi_w$ throughout, 
$t$ denotes a continuous parameter and $w$ a Motzkin word (path). In particular, $\Phi_t$ is a 
is a continuous deformation of $\Phi$. Proceeding as in \eqref{eq:4.7}, we obtain the corresponding moment 
formula for free random variables with respect to the free product of the deformed functionals $\varphi_{i,t}$.
\vspace{3pt}
\begin{Definition}
{\rm By an {\it infinitesimal deformation of Motzkin functionals} we will understand 
the family of multilinear moment functionals $\Phi_{w,t}$ defined by the 
multilinear extension of
{
\setlength{\belowdisplayshortskip}{\baselineskip}
\begin{equation}\label{eq:4.11}
\Phi_{w,t}(a_1, \ldots, a_n)=\Phi_{t}(a_1(j_1)\cdots a_n(j_n)),
\end{equation}
}%
where $w=j_1\cdots j_n\in \mathpzc{M}_n$ and 
$a_1\in \mathcal{A}_{i_1}, \ldots, a_n\in \mathcal{A}_{i_n}$, with 
$n\in \mathbb{N}$. 
}
\end{Definition}
\vspace{3pt}

\begin{Example}
{\rm
Let us give some sample computations of moments with respect to $\Phi_{w,t}$.
If $w=1^3$, we get
\[
\Phi_{w,t}(a_1,a_2,a_3)=\varphi_{i_1,t}(a_1)\varphi_{i_2,t}(a_2)\varphi_{i_2,t}(a_3)
\]
when $i_1\neq i_2\neq i_3$. In turn, if $w=121$, we get
\begin{align*} 
\Phi_{w,t}(a_1,a_2,a_3)&=\widetilde{\varphi}_{i_1,t}(a_1p_{i_1}^{\perp}a_3)\widetilde{\varphi}_{i_2,t}(p_{i_2}a_2p_{i_2})\\
&=\beta_{2,t}(a_1,a_3)\varphi_{i_2,t}(a_2)
\end{align*}
when $i_1=i_3\neq i_2$. Finally, if $w=12121$, we get
\begin{align*} 
\Phi_{w,t}(a_1,a_2,a_3,a_4,a_5)&=\widetilde{\varphi}_{i_1,t}(a_1p_{i_1}^{\perp}a_3p_{i_1}^{\perp}a_5)
\widetilde{\varphi}_{i_2,t}(p_{i_2}a_2p_{i_2}a_4p_{i_2})\\
&=\beta_{3,t}(a_1,a_3,a_5)\varphi_{i_2,t}(a_2)\varphi_{i_2,t}(a_4)
\end{align*}
when $i_1=i_3=i_5\neq i_2=i_4$.
We used \eqref{eq:2.6} (with $j$ replaced by $p_{i_1}$) and \eqref{eq:2.8}. 
In general, an idempotent $p_{i}$ 
factorizes the moment with respect to $\widetilde{\varphi}_{i}$, 
whereas $p_{i}^{\perp}$ transforms its portion into a Boolean cumulant.
}
\end{Example}

It is clear that using the infinitesimal deformations of Motzkin functionals, 
we obtain a decomposition of the free product of deformed functionals. 

\begin{Corollary}
Under the assumptions of Theorem 4.1, it holds that
\begin{equation}\label{eq:4.12}
\varphi_{t}(a_1\cdots a_n)=
\sum_{w=j_1\cdots j_n\in \mathpzc{M}_{n}}
\Phi_{w,t}(a_{1},\ldots ,a_{n}),
\end{equation}
where $\varphi_{t}=(\star_{i\in I} \varphi_{i,t})\circ \tau$, 
whenever $a_1\in \mathcal{A}_{i_1}, \ldots, a_n\in \mathcal{A}_{i_n}$ 
and $i_1\neq \cdots \neq i_n$, $n\in \mathbb{N}$.
\end{Corollary}
{\it Proof.}
This formula follows immediately from Theorem 4.1 when we substitute 
$\varphi_{i,t}$ given by \eqref{eq:4.9} for each $\varphi_i$ in \eqref{eq:4.4}, with the 
use of \eqref{eq:2.8}, where $i\in I$.
\hfill $\blacksquare$\\

\section{Partial tensor evaluations}

We would like to compute $\Phi_{w,t}(a_{1},\ldots ,a_{n})$ and their 
derivatives with respect to $t$ at $t=0$. This will enable us to obtain a new 
approach to infinitesimal freeness and infinitesimal higher-order freeness.

In a natural way, these computations will involve moments and Boolean cumulants 
of orthogonal replicas of variables from the algebras ${\mathcal A}_{i}$. 
Such replicas are indexed by a label $i\in I$ and a color $j\in {\mathbb N}$.
From the definition of the tensor product functional $\Phi_{t}$, 
copies of variables $a_{1}\in {\mathcal A}_{i_1}, \ldots, a_n\in {\mathcal A}_{i_n}$,
denoted $a_{1}(j_1), \ldots, a_{n}(j_n)$,
occupy the same tensor position if and only if both their labels and colors coincide.
Consequently, the computation of a moment of orthogonal replicas depends
on the tuples of labels $(i_1, \ldots, i_n)$ and 
colors $(j_1, \ldots, j_n)$. In particular, if indices $k,m\in [n]$ belong to the same block $V$ 
of a partition $\pi$ associated with a moment with respect to $\Phi_{t}$, it is necessary 
that $j_k=j_m$ and $i_k=i_m$. These conditions, however, are not sufficient since 
idempotents may appear between variables and can separate them even when labels and colors agree.
To encode and control these separations, we study the geometry of Motzkin paths in greater detail.

To put it differently, in order to obtain a nonzero contribution from a given 
mixed moment $\Phi_{w,t}(a_1, \ldots, a_n)$, the level return partition $\pi(w)$ 
associated with the tuples $(i_1, \ldots, i_n)$ and $(j_1, \ldots, j_n)$ must be 
adapted to $(i_1, \ldots, i_n)$. Let us recall that this means that 
$i_k=i_m$ whenever indices $k$ and $m$ belong to the same block $V\in \pi(w)$ and, moreover, 
they satisfy the nested alternation condition (Definition 3.6).
Consequently, we will write the Boolean cumulants associated with a block $V$ 
as $\beta_{|V|,t}((a_{k})_{k\in V})$, where $a_k\in {\mathcal A}_{i}$ for some $i$ 
and for each $k\in V$. In the computations, we will use the formula for Boolean cumulants \eqref{eq:2.6}.

In order to determine the partition $\pi$ of $[n]$ that may 
give a nonzero contribution to the given mixed moment w.r.t. $\Phi_{w,t}$ and then 
compute this contribution, we will examine the level return partition 
$\pi(w)$ associated with $w$. Clearly, it is a subpartition of the partition defined 
by the tuple of colors encoded by $w$, namely,
{
\setlength{\belowdisplayshortskip}{\baselineskip}
\[
S(w)=\{S_{j}: 1\leq j\leq h(w)\}
\]
}%
where $w=j_1\cdots j_n\in \mathpzc{M}_{n}$, $n\in \mathbb{N}$.
We will show that the partition which allows us to compute
$\Phi_{w,t}(a_1, \ldots, a_{n})$ is precisely the level return partition $\pi(w)$ associated 
with $w$. We will assume that neighboring variables have different labels, 
that is $i_1\neq \cdots \neq i_n$, as it is usually assumed in computations of 
moments of noncommutative variables that belong to different algebras.

In order to compute the moments w.r.t. $\Phi_{w,t}$, let us first study
a simple model case, in which $I=\{1,2\}$ 
and we have replicas of just two colors, thus $j_1, \ldots, j_n\in \{1,2\}$ and such that
the associated reduced Motzkin word has just one non-singleton level return block.
\vspace{3pt}
\begin{Lemma}
Let $a_{1}\in {\mathcal A}_{i_1}, \ldots, a_n\in {\mathcal A}_{i_n}$, 
where $i_1\neq\cdots \neq i_n$ and $i_1,\ldots, i_n\in \{1,2\}$,  
and let $w=j_1\cdots j_n\in \mathpzc{M}_{n}$, where $j_1, \ldots, j_n\in \{1,2\}$ be such 
that there is only one non-singleton level return block $V=\{k_1<\cdots <k_q\}$ of color $1$. 
If $i_{k_1}=\cdots =i_{k_q}$, then 
\begin{equation}\label{eq:5.1}
\Phi_{w,t}(a_1, \ldots, a_n)=
\beta_{q,t}(a_{k_1}, \ldots, a_{k_q}) \prod_{{\rm singletons}\;\{k\}}\varphi_{i_k,t}(a_k).
\end{equation}
If not all labels associated with $V$ are identical, 
then $\Phi_{w,t}(a_1, \ldots, a_n)=0$.
\end{Lemma}
{\it Proof.}
Since $w$ is of height $2$, all replicas of color $2$ produce variables 
$a_{k}$ at tensor positions of color $2$ which are separated by idempotents since 
the nearest neighbors of each such replica, $a_{k-1}(j_{k-1}), a_{k+1}(j_{k+1})$ 
have labels different from $i_k$. Therefore, we obtain $p_{i_{k}}a_{k}p_{i_{k}}$ 
at site of color $2$ and label $i_k$. 
If $J_2\subset [n]$ is the set of indices that give singletons of color $2$, we obtain
the nontrivial contribution to $\Phi_{w,t}(a_1, \ldots, a_n)$ of the form
\[
\prod_{k\in J_2}\varphi_{i_{k},t}(a_{k})
\]
when computing the partial tensor product 
$\otimes_{i\in I}\widetilde{\varphi}_{i,t}$ at level $2$.
It remains to compute the remaining tensor product 
$\otimes_{i\in I}\widetilde{\varphi}_{i,t}$ at level $1$.
However, at level $1$ we have all elements from the level return block 
$V=\{k_1<\cdots <k_q\}$ and perhaps, in addition, 
some singletons of color $1$. Let $J_1\subset [n]$ be the set of these singletons. 
As in the case of singletons of color $2$, they give the nontrivial contribution
\[
\prod_{k\in J_1}\varphi_{i_{k},t}(a_{k}).
\]
In turn, the non-singleton level return block $V$ contributes 
\[
\varphi_{i_{k_1},t}(a_{k_1}p_{i_{k_1}}^{\perp}a_{k_2}p_{i_{k_1},t}^{\perp}\cdots p_{i_{k_1},t}^{\perp}a_{k_q})
=\beta_{q,t}(a_{k_1}, \ldots, a_{k_q}), 
\]
provided $i_{k_1}=\cdots =i_{k_q}$. If any two neighboring indices from $V$, say $k_p,k_{p+1}$,
are associated with replicas with different labels, $i_{k_p}\neq i_{k_{p+1}}$, 
then we obtain that $p_{i_{k_p}}^{\perp}$ meets $p_{i_{k_p}}$ which makes 
the moment $\Phi_{w,t}(a_1, \ldots, a_n)=0$ vanish.
Collecting the contributions from all singletons and from the level return block 
gives the desired expression, which completes the proof.
\hfill$\blacksquare$\\

\begin{Example}
{\rm Let $w$ be a pyramid word $w=12\cdots (m-1)m(m+1)\cdots 21$. We already discussed 
its combinatorics in Example 3.1. Assume that we have only two labels, $1$ and $2$ and they alternate in the tuple
$(i_1, \ldots, i_n)$, where $n=2m-1$. 
Using \eqref{eq:5.1}, we obtain
{
\setlength{\belowdisplayshortskip}{\baselineskip}
\begin{align*}
\Phi_{w,t}(a_1\cdots a_n)&=
%\Phi_{t}((a_1p_{i_1}^{\perp}a_n)_{i_1,1}\otimes \cdots \otimes 
%(a_{m-1}p_{i_{m-1}}^{\perp}a_{m+1})_{i_{m-1},m-1}\otimes (a_{m})_{i_m,m})\\
%&= \varphi_{i_1,t}(a_1p_{i_1}^{\perp}a_n)\cdots \varphi_{i_{m-1},t}(a_{m-1}p_{i_{m-1}}^{\perp}a_{m+1})\varphi_{i_{m},t}(a_{m})\\
\beta_{2,t}(a_1, a_n)\cdots \beta_{2,t}(a_{m-1}, a_{m+1})\varphi_{i_m,t}(a_m)
\end{align*}
}%
%where we adopted the convention that only those tensor sites are shown where we have non-trivial variables, with 
%$(a)_{i,j}$ standing for $a$ at the tensor site of label $i$ and color $j$.
Anticipating that we will later attempt to compute the derivative of 
this product of factors w.r.t. $t$ at $t=0$, we can observe that 
if the variables are centered, the only term that survives is the one which applies derivative to $\varphi_{i_m,t}(a_m)$ since $\varphi_{i_m,0}(a_m)=0$. Therefore, 
without any sophisticated computations we can see that
{
\setlength{\belowdisplayshortskip}{\baselineskip}
\begin{align*}
\left.\left[\frac{d\Phi_{w,t}(a_1\cdots a_n)}{dt}\right]\right|_{t=0}&=\beta_{2}(a_1, a_n)\cdots \beta_{2}(a_{m-1}, a_{m+1})\varphi_{i_m}'(a_m)\\
&=\varphi_{i_1}(a_1a_n)\cdots \varphi_{i_{m-1}}(a_{m-1}a_{m+1})\varphi_{i_m}'(a_m)
\end{align*}
}%
which is exactly the formula obtained for the moment of $\varphi'(a_1\cdots a_n)$ in infinitesimal freeness. The reason why 
we obtain this expression from computing $\Phi_{w}'(a_1, \ldots, a_n)$ for just one Motzkin path is that the remaining paths
give zero contribution (this will be showed later).
}
\end{Example}

In order to treat the case of an arbitrary index set, we will first examine
the partial evaluation of the product of replicas associated with a level return block
and the excursions between the level return points.

\begin{Lemma}
Let $a_{1}\in {\mathcal A}_{i_1}, \ldots, a_n\in {\mathcal A}_{i_n}$, where $i_1\neq\cdots \neq i_n$, 
and let $w=j_1\cdots j_n\in \mathpzc{M}_{n}$ be of height $j+1$. If $V=\{k_1<\cdots <k_q\}$ 
is a level return block of color $j$ and $i_{k_1}=\cdots =i_{k_q}$, then the simple tensor
\begin{equation}\label{eq:5.2}
(a_{k_1}p_{i_{k_1}}^{\perp}\cdots p_{i_{k_1}}^{\perp}a_{k_{q}}\otimes p_{i_{k_1}})
\otimes \bigotimes_{\stackrel{k_1<k<k_q}{\scriptscriptstyle k\notin V}}(p_{i_k}\otimes p_{i_k}a_{k}p_{i_k})
\end{equation}
is the partial evaluation of the product of replicas $\prod_{k=k_1}^{k_q}a_{k}(j_k)$ on tensor positions
with labels $i_{k}$, where $k_1\leq k<k_q$, and colors $j,j+1$.
\end{Lemma}
{\it Proof.}
First, consider the product of replicas associated with two level return points $k-1,m+1$ and 
the excursion between them $w'=j_k\ldots j_m$ of the form
\[
a_{k-1}(j_{k-1})[a_{k}(j_k)\cdots a_{m}(j_m)]a_{m+1}(j_{m+1})
\]
where $j_{k-1}=j_{m+1}=j$ and $j_k=\cdots = j_m=j+1$. By assumption, $i_{k-1}=i_{m-1}\neq i_{s}$
for all $k\leq s\leq m$. Suppose first that all labels $i_s$ in the excursion are different from each other. 
Without loss of generality we can assume that they are 
increasing, that is $i_k<i_{k+1}<\cdots <i_m$. Otherwise we can rearrange the order in which 
we write tensor products. 
Let us carry out the partial evaluation of the above product of replicas 
(treated as a tensor product) on tensor positions of colors $j-1,j$ and 
labels $i_{k-1}=i_{m+1},i_k,\ldots, i_m$. Computing this partial evaluation, we obtain the product
\[
\prod_{p=k-1}^{m+1}T_{p}\in \bigotimes_{p=k-1}^{m}
(\widetilde{\mathcal A}_{i_p}\otimes \widetilde{\mathcal A}_{i_p})
\]
where 
\begin{align*}
T_{k-1}&=(a_{k-1}\otimes 1_{i_{k-1}})\,\otimes (p_{i_k}\otimes p_{i_k})
\,\otimes \cdots \otimes (p_{i_m}\otimes p_{i_m})\\
T_{k}&=(p_{i_{k-1}}^{\sigma_{1,0}}\otimes p_{i_{k-1}})\otimes (1_{i_k}\otimes a_{k})
\;\,\otimes \cdots \otimes (p_{i_m}^{\sigma_{1,r}}\!\!
\otimes p_{i_m})\\
T_{k+1}&=(p_{i_{k-1}}^{\sigma_{2,0}}\otimes p_{i_{k-1}})\otimes (p_{i_k}^{\sigma_{2,1}}\!\!\otimes p_{i_k})
\otimes \cdots \otimes (p_{i_m}^{\sigma_{2,r}}\!\!
\otimes p_{i_m})\\
&\vdots\\
T_{m}&= (p_{i_{k-1}}^{\sigma_{r,0}}\otimes p_{i_{k-1}})\otimes (p_{i_k}^{\sigma_{r,1}}\!\!\otimes p_{i_k})
\otimes \cdots \otimes (1_{i_m}
\otimes a_{m})\\
T_{m+1}&= (a_{m+1}\otimes 1_{i_{k-1}})\otimes (p_{i_k}\otimes p_{i_k})
\,\otimes \cdots \otimes (p_{i_{m}}
\otimes p_{i_m})
\end{align*}
where $r=m-k+1$ and $\sigma_{s,t}\in \{1, \perp\}$ for all $s,t$ and in each row 
there is at least one $\perp$. The last condition follows from the formula
\[
\bigotimes_{r\in J}1_{i_{r}}-\bigotimes_{r\in J}p_{i_{r}}
=\sum_{\stackrel{\sigma_{k}\in \{1,\perp\}}{\scriptscriptstyle k\in J}}^{\circ}\bigotimes_{r\in J}p_{i_r}^{\sigma_{r}}.
\]
for any finite set $J$, where the circle above the sum means that at least one $\sigma_{k}$ is equal to $\perp$.
In order to get a nonzero contribution we have to rule out the configurations in which $p_{i}$ 
meets $p_{i}^{\perp}$ for any $i$. When we multiply projections along all `columns' but the first one, we get 
the set of conditions
\[
\sigma_{2,1}=\sigma_{3,1}=\cdots = \sigma_{r,1}\;\;=1
\]
\[
\sigma_{1,2}=\sigma_{3,2}=\cdots = \sigma_{r,2}\;\;=1
\]
\[
\vdots
\]
\[
\sigma_{1,r}=\sigma_{2,r}=\cdots = \sigma_{r,r-1}=1
\]
which implies that 
\vspace{5pt}
\[
\sigma_{1,0}=\sigma_{2,0}=\cdots =\sigma_{r,0}=\perp
\]
\vspace{5pt}
since in each row we must have at least one $\perp$. Therefore, we obtain
\[
\prod_{p=k-1}^{m+1}T_{p}=
(a_{k-1}p_{i_{k-1}}^{\perp}a_{m+1}\otimes p_{i_{k-1}})\otimes (p_{i_k}\otimes p_{i_k}a_{k}p_{i_k})
\otimes \cdots \otimes (p_{i_m}\otimes p_{i_m}a_{m}p_{i_m}).
\]
It is not hard to see that if some labels in an excursion between 
$k$ and $m$ are identical then we obtain the same expression. 
Repeating this evaluation for all elements of a level return block 
$V=\{k_1< \cdots <k_q\}$ and the excursions between them, we obtain 
\[
\prod_{p=k_1}^{k_q}T_{p}=
(a_{k_1}p_{i_{k_1}}^{\perp}\cdots p_{i_{k_1}}^{\perp}a_{k_{q}}\otimes p_{i_{k_1}})
\otimes \bigotimes_{\stackrel{k_1<k<k_q}{\scriptscriptstyle k\notin V}}(p_{i_k}\otimes p_{i_k}a_{k}p_{i_k}).
\]
This completes the proof of \eqref{eq:5.2}. 
\hfill $\blacksquare$\\

\begin{Lemma}
Let $a_{1}\in {\mathcal A}_{i_1}, \ldots, a_n\in {\mathcal A}_{i_n}$, where $i_1\neq\cdots \neq i_n$, 
and let $w=j_1\cdots j_n\in \mathpzc{M}_{n}$ be of height $j+1$. If $V=\{k_1<\cdots <k_q\}$ 
is a level return block of color $j$, then 
\begin{equation}\label{eq:5.3}
\prod_{k=k_1}^{k_q}a_{k}(j_k)=0
\end{equation}
when not all associated labels $i_{k_1},\ldots, i_{k_q}$ are equal to each other.
\end{Lemma}
{\it Proof.} 
Let us suppose that two consecutive level return points associated with different lables 
are $k-1, m+1$. Moreover, assume again (without loss of generality) 
that the labels in the excursion are increasing, that is $i_k<i_{k+1}<\cdots <i_m$.
Using similar notations to those in the proof of Lemma 5.2, 
the partial evaluation of the product of replicas
\[
a_{k-1}(j_{k-1})[a_{k}(j_k)\cdots a_{m}(j_m)]a_{m+1}(j_{m+1})
\]
gives in this case
\[
\prod_{p=k-1}^{m+1}T_{p}\in \bigotimes_{p=k-1}^{m+1}
(\widetilde{\mathcal A}_{i_p}\otimes \widetilde{\mathcal A}_{i_p})
\]
where 
{
\setlength{\belowdisplayshortskip}{\baselineskip}
\begin{align*}
T_{k-1}&=(a_{k-1}\otimes 1_{i_{k-1}})\,\otimes (p_{i_k}\otimes p_{i_k})\,
\otimes \cdots \otimes (p_{i_m}\otimes p_{i_m})\otimes (p_{i_{m+1}}\otimes p_{i_{m+1}})\\
T_{k}&=(p_{i_{k-1}}^{\sigma_{1,0}}\otimes p_{i_{k-1}}) \otimes (1_{i_k}\otimes a_{k})
\otimes \cdots \otimes (p_{i_m}^{\sigma_{1,r}}\!\!\otimes p_{i_m})\otimes 
(p_{i_{m+1}}^{\sigma_{1,r+1}}\!\!\otimes p_{i_{m+1}})\\
&\vdots\\
T_{m}&= (p_{i_{k-1}}^{\sigma_{r,0}}\otimes p_{i_{k-1}})\otimes (p_{i_k}^{\sigma_{r,1}}\!\!\otimes p_{i_k})
\otimes \cdots \otimes (1_{i_m}
\otimes a_{m})\otimes (p_{i_{m+1}}^{\sigma_{r,r+1}}\otimes p_{i_{m+1}})\\
T_{m+1}&= (p_{i_{k-1}}\otimes p_{i_{k-1}}) \otimes (p_{i_k}\otimes p_{i_k})
\otimes \cdots \otimes (p_{i_{m}} \otimes p_{i_m})\otimes (a_{m+1}\otimes 1_{i_{m+1}}).
\end{align*}
}%
Note that in this case, in contrast to the situation of Lemma 5.2, 
variables $a_{k-1}$ and $a_{m+1}$ occupy different tensor positions
since $i_{k-1}\neq i_{m+1}$. Therefore, we obtain one more column, 
and thus, in addition to the conditions on $\sigma_{i,j}$ given in the proof of Lemma 5.2,
we have: 
\[
\sigma_{1,r+1}=\sigma_{2,r+1}=\cdots = \sigma_{r,r+1}=1.
\]
Again, this leads to
\[
\sigma_{1,0}=\cdots = \sigma_{r,0}=\perp
\] 
since in each row we must have at least one $\perp$. However,
in the present situation this means that $p_{i_{k-1}}$ meets $p^{\perp}_{i_{k-1}}$ 
in the first column since we have $p_{i_{k-1}}$ in the lowest left corner
instead of $a_{m+1}$.
Therefore, the product of replicas vanishes, which proves \eqref{eq:5.3}.
\hfill$\blacksquare$\\

\section{Formula for moments} 

Partial tensor evaluations of Lemmas 5.2-5.3 give the foundations for computations of moments
$\Phi_{w,t}(a_1, \ldots, a_n)$. 
These evaluations were localized by the level return blocks and the excursions between their elements.
Moreover, we assumed in Lemma 5.2 that these blocks were of color 
$j$, where $j+1$ was the height of $w$, but it is not hard to see that
similar evaluations can be performed at all levels.

Since we begin with the highest colors present in 
$w$ and then proceed downward to color $1$, one can state the results 
using tensor product conditional expectations.
For a sequence $(B_n)_{n\in {\mathbb N}}$ of unital algebras with
associated sequence of normalized linear functionals $(\psi_{n})_{n\in {\mathbb N}}$, 
consider the tensor product 
\[
B(m):=B_1\otimes \cdots \otimes B_m,\;\;\;m\in {\mathbb N}
\] 
with the convention $B(0)={\mathbb C}$. Let
\[
\Psi_{m}:=\psi_1\otimes \cdots \otimes \psi_m, \;\;\;m\in {\mathbb N},
\]
denote the corresponding products functionals on $B(m)$.

We define a {\it tower of conditional expectations} to be the sequence of linear maps
$E_{m}:B(m)\rightarrow B(m-1)$, where $m\in {\mathbb N}$, given on simple tensors by
\begin{align*}
E_{m}(b_1\otimes \cdots \otimes b_m)&=\psi_m(b_m)(b_1\otimes \cdots \otimes b_{m-1}).
\end{align*}
Iterating these maps gives 
\[
\Psi_m(b_1\otimes \cdots \otimes b_m)=E_1(E_2(\ldots E_m(b_1\otimes \cdots \otimes b_m)\ldots )).
\]

In our setting, the algebras $B_j$ are themselves tensor products, namely
\[
B_j=\bigotimes_{i\in I}\widetilde{\mathcal A}_{i}\;\;\;{\rm and}\;\;\;
\psi_j=\bigotimes_{i\in I}\widetilde{\varphi}_{i,t}
\] 
for any $j\in {\mathbb N}$ and we interpret the index $j$ as a color. When 
computing moments $\Phi_{w,t}(a_1, \ldots, a_m)$ for $w=j_1\cdots j_n\in \mathpzc{M}_{n}$
with the use of the tower of conditional expectations, the computation 
reduces to the algebra
\[
B(j)=\bigotimes_{1\leq k \leq j}B_k\cong \bigotimes_{1\leq k\leq j}\bigotimes_{i\in I}\widetilde{A}_{i}
\cong \bigotimes_{i\in I}\widetilde{A}_{i}^{\otimes j}
\]
where $j=h(w)$ is the height of $w$. The reduction is possible because at all remaining tensor positions
we insert idempotents $p_i\in \widetilde{\mathcal{A}}_{i}$ for which we assume 
$\widetilde{\varphi}_{i,t}(p_i)=1$, 
so these factors do not affect the value of the conditional expectations.

For a finite set $V=\{k_1<\ldots <k_r\}$, the family $(a_{k})_{k\in V}$ 
will always be understood as the ordered tuple $(a_{k_1}, \ldots, a_{k_r})$.
This notation will be used in the context of various cumulants. 
For instance, Boolean cumulants will be denoted by $\beta_{|V|}((a_{k})_{k\in V})$, which is equivalent
to some other authors' $\beta_{|V|}(a_1, \ldots, a_n|V)$ or $\beta_{|V|}(a_{V})$.

\begin{Theorem}
Let $a_1\in {\mathcal A}_{i_1}, \ldots, a_n\in \mathcal{A}_{i_n}$, where $i_1\neq \cdots \neq i_n$
and let $w\in \mathpzc{M}_{n}$.
\begin{enumerate}
\item
If $w$ is adapted to $(i_1, \ldots, i_n)$, then
\begin{equation}\label{eq:6.1}
\Phi_{w,t}(a_1, \ldots, a_n)=
\prod_{V\in \pi(w)}\beta_{|V|,t}((a_{i})_{i\in V}),
\end{equation}
where $\pi(w)$ is the level return partition associated with $w$.
\item
If $w$ is not adapted to $(i_1, \ldots, i_n)$, then $\Phi_{w,t}(a_1, \ldots, a_n)=0$.
\end{enumerate}
\end{Theorem}

{\it Proof.}
Let $w=j_1\ldots j_n$ be of height $j+1$, as in the proof of Lemma 5.2. 
The partial tensor evaluation of that Lemma tells us that if $j_k=j+1$, 
then the corresponding replica contributes 
$p_{i_k}a_{k}p_{i_k}$. This holds even if $k=1$ or $k=n$ since we are free to write 
the idempotent $p_{i_k}$ before $a_{1}$ and after $a_n$ since 
$\widetilde{\varphi}_{i_{k}}(p_{i_k}up_{i_k})=\widetilde{\varphi}_{i_{k}}(u)$ for any
$u\in \widetilde{{\mathcal A}}_{i_k}$. Besides, 
{
\setlength{\belowdisplayshortskip}{\baselineskip}
\[
\widetilde{\varphi}_{i_k,t}(u_1pa_kpu_2)=
\varphi_{i_k,t}(a_k)\widetilde{\varphi}_{i_k,t}(u_1)
\widetilde{\varphi}_{i_k,t}(u_2)
\]
}%
\noindent
by \eqref{eq:2.8}. It is worth pointing out that each such singleton 
$\{k\}$ gives a local maximum of $w$. Namely, we have
$j_{k-1}\leq j_k\geq j_{k+1}$, with the understanding that if $k=1$ or $k=n$ only one inequality 
holds. Taking into account all local maxima at level $j+1$ and computing 
the conditional expectation $E_{j+1}$, we obtain 
\[
E_{j+1}[a_{1}(j_1)\cdots a_{n}(j_n)]_{j+1}=
\Bigl\lgroup\prod_{\stackrel{{\rm singleton\; blocks}\; \{k\}}{\scriptscriptstyle {\rm of\;color}\;j+1}}
{\varphi}_{i_k}(a_{k})\Bigr\rgroup[a_{1}(j_1)\cdots a_{n}(j_n)]_{j}
\]
where $[a]_{h}$ is the truncation of $a$ to the tensor product of order $h$.
Now, at level $j$ we have level return blocks of color $j$ (they can be singletons).
Using again Lemma 5.2, we obtain
\[ 
E_{j}[(a_{1}(j_1)\cdots a_{n}(j_n))]_{j}=
\Bigl\lgroup\prod_{\stackrel{{\rm level\;return\;blocks}\;V}{\scriptscriptstyle{\rm of\;color}\;j\;{\rm adapted\;to\;labels} }}
\beta_{|V|,t}((a_{i})_{i\in V})\Bigr\rgroup
[(a_{1}(j_1)\cdots a_{n}(j_n))]_{j-1}
\]
where, in particular, some of these level return blocks can be singleton blocks. 
We continue in this fashion with the remaining conditional expectations.
After applying the whole tower of conditional expectations 
$(E_{k})_{1\leq k \leq j+1}$ we obtain the product of 
Boolean cumulants corresponding to all level return blocks at all levels. 
In particular, if $\{k\}$ is a singleton block, then $\beta_{1,t}(a_{k})=
\varphi_{i_k}(a_{k})$. This gives the desired formula \eqref{eq:6.1} since
\vspace{3pt}
{
\setlength{\belowdisplayshortskip}{\baselineskip}
\begin{align*}
\Phi_{w,t}(a_1, \ldots, a_n)&=
E_1[E_2[\ldots [E_{j+1}[a_{1}(j_1)\cdots a_{n}(j_n)]_{j+1}]\ldots ]_2]_1\\
&=
\prod_{V\in \pi(w)}\beta_{|V|,t}((a_{i})_{i\in V})
\end{align*}
}%
whenever $w$ is adapted to the tuple $\ell=(i_1, \ldots, i_n)$, 
which follows from Lemma 5.2. Now, if $w$ is not adapted to 
$(i_1, \ldots, i_n)$, then there exists a level return block 
$V=\{k_1<\cdots <k_q\}\in \pi(w)$ for which the labels $i_{k_1}, \ldots, i_{k_{q}}$ 
are not identical. In that case, we get zero conditional expectation at the highest 
level at which such a situation takes place by Lemma 5.3. Therefore, 
in that case $\Phi_{w,t}(a_{1}, \ldots, a_{n})=0$, which gives (2). This completes the proof.
\hfill $\blacksquare$\\

\begin{Corollary}
The formula of Theorem 6.1 can be written in the equivalent form
\begin{equation}\label{eq:6.2}
\Phi_{w,t}(a_1, \ldots, a_n)=\prod_{{\rm singleton\;blocks} \;\{k\}}\varphi_{i_k,t}(a_{k})
\prod_{{\rm non-singleton\; blocks\;V}}\beta_{|V|,t}((a_{i})_{i\in V})
\end{equation}
whenever $w$ is adapted to $(i_1, \ldots, i_n)$. 
\end{Corollary}
{\it Proof.}
This formula follows immediately from \eqref{eq:6.1} if we distinguish singletons from the remaining 
level return blocks.
\hfill $\blacksquare$\\

\section{Infinitesimal Motzkin functionals} 

We are going to compute the derivatives of functionals $\Phi_{w,t}$ at $t=0$. It is interesting to 
note that the only functional which has a non-zero derivative is that associated with a pyramid Motzkin path.

We begin with the definition of $\Phi_{w}'$ and the definition of the pyramid Motzkin path.

\begin{Definition}
{\rm By {\it infinitesimal Motzkin functionals} we will understand the family of 
multilinear moment functionals defined by the multilinear extension of
\vspace{3pt}
{
\setlength{\belowdisplayshortskip}{\baselineskip}
\begin{equation}\label{eq:7.1}
\Phi_{w}'\;(a_1, \ldots, a_n)=\left.\left[\frac{d\Phi_{w,t}(a_1, \ldots, a_n)}{dt}\right]\right|_{t=0}
\end{equation}
}%
where $w=j_1\cdots j_n\in \mathpzc{M}_{n}$ and 
$a_1\in {\mathcal A}_{i_1}, \ldots , a_n\in {\mathcal A}_{i_n}$, with $n\in {\mathbb N}$.
}
\end{Definition}

\begin{Definition}
{\rm A Motzkin path is called a {\it pyramid Motzkin path of length} $2m$ (or, simply, 
a {\it pyramid path}) if and only if its step word is of the form 
\[
U^{m}D^{m},
\]
where $m\geq 1$. The corrresponding $w\in \mathpzc{M}$, called a 
{\it pyramid word of length $2m+1$}, is of the form 
$w=1\cdots m(m+1)m\cdots 1$. The empty path (corresponding to the word $1$) 
will be treated as the Motzkin pyramid path of length $0$.}
\end{Definition}

\begin{Lemma}
A Motzkin path has exactly one local maximum if and only if it is of even length and 
it is a pyramid Motzkin path (the corresponding $w\in \mathpzc{M}_{n}$ is of odd length). 
\end{Lemma}
{\it Proof.}
The case of the empty Motzkin path is formally treated as a path with 
one local maximum. A non-empty Motzkin path is a pyramid Motzkin path if and only if it is of even length 
(the corresponding word $w$ is of odd length). Then is has exactly one local maximum. 
If it is not a pyramid Motzkin path, then it has at least two local maxima since 
we have three cases:
\begin{enumerate}
\item it has at least two peaks, 
\item it has at least one plateau of the form 
\[
w'=j_{k}\cdots j_{l}=j(j+1)\cdots (j+1)j
\] 
(with local maxima at $(k+1,j+1)$ and at $(l-1,j+1)$),
\item
it has at least one peak and at least one plateau,
\[
w'=j_{k}\cdots j_l=j(j+1)\cdots (j+1)(j+2)
\] 
(with local maximum at $(k+1,j+1)$), or
\[
w'=j_{k}\cdots j_l= (j+2)(j+1)\cdots (j+1)j
\] 
(with local maximum at $(l-1,j+1)$).
\end{enumerate}
This completes the proof.
\hfill $\blacksquare$\\

Let us first compute the infinitesimal Motzkin functionals for pyramid Motzkin paths in the 
case when the variables are centered.

\begin{Lemma}
Let $a_1\in {\mathcal A}_{i_1}^{\circ}, \ldots, a_n\in \mathcal{A}_{i_n}^{\circ}$, 
where $i_1\neq \cdots \neq i_n$ and let $w\in \mathpzc{M}_{n}$, $n\in {\mathbb N}$.
\begin{enumerate}
\item
If $w$ is a pyramid word of the form $w=1\cdots (m-1)m(m+1)\cdots 1$ and it is adapted to
$(i_1, \ldots, i_n)$, then 
\begin{equation}\label{eq:7.2}
\Phi_{w}'(a_1, \ldots, a_n)=
\varphi_{i_1}(a_1a_{n})\cdots \varphi_{i_{m-1}}(a_{m-1}a_{m+1})\varphi_{i_{m}}'(a_{m}).
\end{equation}
\item
In the remaining cases (if $w$ is not a pyramid Motzkin path or it is not adapted to $(i_1, \ldots, i_n)$),
$\Phi_{w}'(a_1, \ldots, a_n)=0$.
\end{enumerate}
\end{Lemma}
{\it Proof.}
We apply Definition given by \eqref{eq:7.1} to \eqref{eq:6.2}.
If $n=1$, then the statement is obviously true. Therefore, let $n>1$. For each moment associated with a singleton $\{k\}$, we have
$$\varphi_{i_k,0}(a_k)=\varphi_{i_k}(a_k)=0,$$
which means that the path $w$ can have at most one local maximum if the derivative 
$\Phi_w$ is to be non-zero.  The only such Motzkin path
that has only one maximum is the pyramid path $w=1\cdots (m-1)m(m+1)\cdots 1$. 
In this case, we have 
\[
\Phi_{w}'(a_1, \ldots, a_n)=
\varphi_{i_1}(a_1a_{n})\cdots \varphi_{i_{m-1}}(a_{m-1}a_{m+1})\varphi_{i_m}'(a_m)
\]
by Corollary 6.1 whenever $i_1=i_n, \ldots, i_{m-1}=i_{m+1}$ since 
\[
\beta_{2}(a_1,a_n)=\varphi_{i_1}(a_1a_n), \ldots, 
\beta_2(a_{m-1},a_{m+1})=\varphi_{i_{m-1}}(a_{m-1}a_{m+1}).
\]
If $w$ is a pyramid path, but 
it is not adapted to labels, the derivative vanishes. This completes the proof.
\hfill $\blacksquare$\\

\begin{figure}
\unitlength=1mm
\special{em:linewidth 1pt}
\linethickness{0.5pt}
\begin{picture}(180.00,40.00)(-17.00,50.00)

%%%%%%%%%%%%%%%path %%%%%%%%%%%%%%%
\put(-5.00,70.00){$w$}
\put(61.00,70.00){$\pi(w)$}
%%%%%%%%%%%%%%%%% path 1%%%%%%%%%%%%%%%%%%%%%%%%
\put(-5.00,55.00){\line(1,1){5.00}}
\put(0.00,60.00){\line(1,1){5.00}}
\put(5.00,65.00){\line(1,1){5.00}}
\put(10.00,70.00){\line(1,1){5.00}}
\put(15.00,75.00){\line(1,1){5.00}}
\put(20.00,80.00){\line(1,-1){5.00}}
\put(25.00,75.00){\line(1,-1){5.00}}
\put(30.00,70.00){\line(1,-1){5.00}}
\put(35.00,65.00){\line(1,-1){5.00}}
\put(40.00,60.00){\line(1,-1){5.00}}

\put(-5.00,55.00){\circle*{1.00}}
\put(0.00,60.00){\circle*{1.00}}
\put(5.00,65.00){\circle*{1.00}}
\put(10.00,70.00){\circle*{1.00}}
\put(15.00,75.00){\circle*{1.00}}
\put(20.00,80.00){\circle*{1.00}}
\put(20.00,80.00){\circle{3.00}}
\put(25.00,75.00){\circle*{1.00}}
\put(30.00,70.00){\circle*{1.00}}
\put(35.00,65.00){\circle*{1.00}}
\put(40.00,60.00){\circle*{1.00}}
\put(45.00,55.00){\circle*{1.00}}

%%%% partition  %%%%%%%%%%%%%%%%%%%%%%

\put(75.00,55.00){\circle*{1.00}}
\put(80.00,55.00){\circle*{1.00}}
\put(85.00,55.00){\circle*{1.00}}
\put(90.00,55.00){\circle*{1.00}}
\put(95.00,55.00){\circle*{1.00}}
\put(100.00,55.00){\circle*{1.00}}
\put(100.00,55.00){\circle{3.00}}
\put(105.00,55.00){\circle*{1.00}}
\put(110.00,55.00){\circle*{1.00}}
\put(115.00,55.00){\circle*{1.00}}
\put(120.00,55.00){\circle*{1.00}}
\put(125.00,55.00){\circle*{1.00}}

\put(74.00,51.00){$\scriptstyle{1}$}
\put(79.00,51.00){$\scriptstyle{2}$}
\put(84.00,51.00){$\scriptstyle{3}$}
\put(89.00,51.00){$\scriptstyle{4}$}
\put(94.00,51.00){$\scriptstyle{5}$}
\put(99.00,51.00){$\scriptstyle{6}$}
\put(104.00,51.00){$\scriptstyle{7}$}
\put(109.00,51.00){$\scriptstyle{8}$}
\put(114.00,51.00){$\scriptstyle{9}$}
\put(118.50,51.00){$\scriptstyle{10}$}
\put(123.50,51.00){$\scriptstyle{11}$}

\put(75.00,55.00){\line(0,1){25}}
\put(80.00,55.00){\line(0,1){20}}
\put(85.00,55.00){\line(0,1){15}}
\put(90.00,55.00){\line(0,1){10}}
\put(95.00,55.00){\line(0,1){5}}
\put(105.00,55.00){\line(0,1){5}}
\put(110.00,55.00){\line(0,1){10}}
\put(115.00,55.00){\line(0,1){15}}
\put(120.00,55.00){\line(0,1){20}}
\put(125.00,55.00){\line(0,1){25}}

\put(75.00,80.00){\line(1,0){50}}
\put(80.00,75.00){\line(1,0){40}}
\put(85.00,70.00){\line(1,0){30}}
\put(90.00,65.00){\line(1,0){20}}
\put(95.00,60.00){\line(1,0){10}}

\end{picture}
\caption{Pyramid Motzkin path of length $10$ and the corresponding fully nested level return partition
with a central singleton block.}
\end{figure}
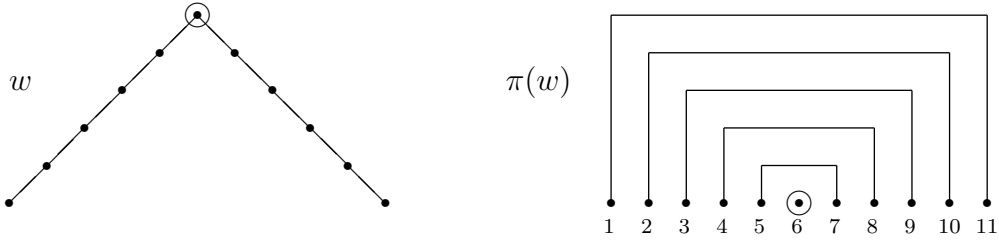

A formula for arbitrary (non-centered) variables is given below. 

\begin{Proposition}
Let $a_1\in {\mathcal A}_{i_1}, \ldots, a_n\in \mathcal{A}_{i_n}$, where $i_1\neq \cdots \neq i_n$
and let $w\in \mathpzc{M}_{n}$. Then it holds that
{
\setlength{\belowdisplayshortskip}{\baselineskip}
\begin{equation}\label{eq:7.3}
\Phi_{w}'(a_1, \ldots, a_n)=\sum_{V\in \pi(w)}\beta_{|V|}'((a_{i})_{i\in V})
\prod_{U\neq V}\beta_{|U|}((a_{i})_{i\in U})
\end{equation}
}%
whenever $w$ is adapted to $(i_1, \ldots, i_n)$ and otherwise it vanishes.
\end{Proposition}
{\it Proof.}
This formula follows immediately from computing the derivative with respect to $t$ of the formula for 
$\Phi_{w,t}$ given by \eqref{eq:6.1} and setting $t=0$.\hfill $\blacksquare$\\

\begin{Theorem}
Let $({\mathcal A}_{i},\varphi_i)_{i\in I}$ be noncommutative probability spaces. 
Let $\varphi'=d\varphi_{t}/dt|_{t=0}$, where $\varphi_{t}=(\star_{i\in I} \varphi_{i,t})\circ \tau$.
Let $a_1\in \mathcal{A}_{i_1}, \ldots, a_n\in \mathcal{A}_{i_n}$, where $i_1\neq \cdots \neq i_n$, $n\in {\mathbb N}$.
\begin{enumerate}
\item 
The following decomposition holds:
\begin{equation}\label{eq:7.4}
\varphi'(a_1\cdots a_n)=
\sum_{w\in \mathpzc{M}_{n}}
\Phi_{w}'(a_{1}, \ldots, a_n).
\end{equation}
\item
If $a_1\in {\mathcal A}_{i_1}^{\circ}, \ldots, a_n\in \mathcal{A}_{i_n}^{\circ}$,  then 
\begin{equation}\label{eq:7.5}
\varphi'(a_1 \cdots a_n)=\varphi_{i_{m}}'(a_{m})\varphi_{i_1}(a_1a_n)\cdots \varphi_{i_m}(a_{m-1}a_{m+1})
\end{equation}
whenever $n=2m+1$ and $i_1=i_n, \ldots, i_{m-1}=i_{m+1}$, and otherwise it vanishes.
\end{enumerate}
\end{Theorem}
{\it Proof.}
The decomposition in (1) follows from Corollary 4.1. 
To obtain (2), it suffices to use (1) and Lemma 7.2. This completes the proof.
\hfill $\blacksquare$

\begin{Example}
{\rm 
Let us consider Motzkin paths of length $4$ (words of length $5$). Let 
$a_k\in {\rm Ker}(\varphi_{i_k})$ for each $k=1, \ldots,5$, where $i_1\neq \cdots \neq i_5$.  
If $w=12321$ (pyramid path), then
{
\setlength{\belowdisplayshortskip}{\baselineskip}
\[
\Phi_{w}'(a_1,a_2,a_3,a_4,a_5)=\varphi_{i_2}(a_1a_5)\varphi_{i_3}(a_2a_4)\varphi_{i_3}'(a_3)
\]
}%
whenever $i_1=i_5$ and $i_2=i_4$ and otherwise it vanishes, since $\beta_{2}(a_1,a_5)=\varphi_{i_1}(a_1a_5)$ 
and $\beta_{2}(a_2,a_4)=\varphi_{i_2}(a_2a_4)$.
The remaining Motzkin paths for $n=5$ give zero contribution since they have more than one local maximum. 
Therefore,
{
\setlength{\belowdisplayshortskip}{\baselineskip}
\[
\varphi'(a_1a_2a_3a_4a_5)=
\delta_{i_1,i_5}\delta_{i_2,i_4}\varphi_{i_2}(a_1a_5)\varphi_{i_3}(a_2a_4)\varphi_{i_3}'(a_3),
\]
}%
where $\varphi'$ is the derivative of the free product $\star_{i\in I}\varphi_{i}$.
}
\end{Example}

We can also consider infinitesimal Boolean independence since the Boolean
independence is included in the framework of Motzkin functionals \cite{L4}.
In that case, it suffices to consider constant Motzkin words $w=1^{n}$ and the associated functionals $\Phi_{1^{n}}$
since
\[
\Phi_{1^n}(a_1, \ldots, a_n)=\varphi_{i_1}(a_1)\cdots \varphi_{i_n}(a_n)
\]
for any $a_1\in {\mathcal A}_{i_1}, \ldots, a_{n}\in {\mathcal A}_{i_n}$, where $i_1\neq \cdots \neq i_n$. 
Using $\pi(w)$ introduced in the context of infinitesimal moments,
all blocks of $\pi(w)$ are singletons if $w=1^{n}$. 
It is therefore natural to expect that if we take the pair $(\Phi_{1^{n}}, \Phi_{1^{n}}')$, then
we will obtain  infinitesimal Boolean independence in our framework. 
Let us remark that the non-centered general version given by Proposition 7.1 is more appropriate since 
in Boolean independence, apart from the fact that units are not identified, 
the factorization formula given above has to hold for any elements, not only the centered ones.

\begin{Definition}
{\rm A Motzkin path is called a {\it flat Motzkin path of length} $n$ (or, simply, 
a {\it flat Motzkin path}) if and only if its step word is of the form $H^{n}$
where $n\geq 0$. The corrresponding reduced Motzkin word, called a 
{\it constant Motzkin word of length $n+1$}, is of the form $w=1^{n+1}$. We allow 
flat empty paths, then $n=0$ and the corresponding word is $w=1$.}
\end{Definition}
\vspace{3pt}

\begin{Definition}
{\rm Let $({\mathcal A}, \varphi, \varphi')$ be an infinitesimal noncommutative probability space.
We say that the family of subalgebras $({\mathcal A}_{i})_{i\in I}$ is {\it infinitesimally Boolean independent}
with respect to $(\varphi, \varphi')$ if and only if 
\begin{enumerate}
\item for any $a_1\in {\mathcal A}_{i_1}, \ldots, a_n\in {\mathcal A}_{i_n}$, where $i_1\neq \ldots \neq i_n$,
it holds that
\begin{equation}\label{eq:7.6}
\varphi'(a_1, \ldots, a_n)=\sum_{k=1}^{n}\varphi(a_{1})\cdots \varphi'(a_k)\cdots\varphi(a_{n})
\end{equation}
\item
the family $({\mathcal A}_{i})_{i\in I}$ is Boolean independent with respect to $\varphi$.
\end{enumerate}
}
\end{Definition}
\vspace{3pt}

In order to obtain infinitesimal Boolean independence in our framework, 
we will use the linear functionals $\chi,\chi'$ on $*_{i\in I}\mathcal A_i$ defined by
\[
\chi(a_1\cdots a_n)=\Phi_{1^n}(a_1,\dots,a_n),
\qquad
\chi'(a_1\cdots a_n)=\Phi_{1^n}'(a_1,\dots,a_n).
\]
whenever $a_1 \in \mathcal A_{i_1}, \ldots, a_n\in {\mathcal A}_{n}$ and $i_1 \neq \cdots \neq i_n$.
Let us recall that $*_{i\in I}\mathcal A_i$ stands 
for the free product of algebras without identification of units.\\

\begin{Proposition}
Let $({\mathcal A}_{i})_{i\in I}$ be a family of unital algebras and 
let $a_1\in {\mathcal A}_{i_1}, \ldots, a_n\in \mathcal{A}_{i_n}$, 
where $i_1\neq \cdots \neq i_n$. Then
\begin{equation}\label{eq:7.7}
\Phi'_{1^n}(a_1,\dots,a_n)
=
\sum_{k=1}^n
\varphi_{i_1}(a_1)\cdots
\varphi_{i_k}'(a_k)\cdots
\varphi_{i_n}(a_n).
\end{equation}
Thus, the family $\left(\mathcal A_i\right)_{i\in I}$ is infinitesimally Boolean independent with respect to 
$(\chi,\chi')$.
\end{Proposition}
{\it Proof.}
For any $a_1\in {\mathcal A}_{i_1}, \ldots , {\mathcal A}_{i_n}$, where $i_1\neq \cdots \neq i_n$, 
we have 
\[
\Phi_{1^{n},t}(a_1, \ldots, a_n)=\varphi_{i_1,t}(a_1)\cdots \varphi_{i_{n},t}(a_{n})
\]
by Corollary 6.1. Therefore,
\[
\Phi_{1^n}'(a_1, \ldots, a_n)=\sum_{k=1}^{n}\varphi_{i_1}(a_1)\cdots \varphi_{i_k}'(a_k)\cdots \varphi_{i_n}(a_n)
\]
which proves our assertion. 
\hfill $\blacksquare$\\

Interestingly enough, we can also write an equivalent formula of `one local maximum type'. 
If we take each variable to be either centered or an internal unit, then the infinitesimal 
Motzkin functional associated with a flat path vanishes unless
there is exactly one centered variable and the remaining ones are internal units. 

\begin{Corollary}
Let $({\mathcal A}_{i})_{i\in I}$ be a family of unital algebras and 
let $a_1\in {\mathcal A}_{i_1}^{\circ}\cup \{1_{i_1}\}, \ldots, a_n\in {\mathcal A}_{i_n}^{\circ}\cup \{1_{i_n}\}$, where
$i_1\neq \cdots \neq i_n$. Then
\begin{equation}\label{eq:7.8}
\Phi_{1^n}'(a_1, \ldots, a_n)=\varphi_{i_{k}}'(a_{k})
\end{equation}
whenever there exists exactly one $k$ for which $a_k\in {\mathcal A}_{i_k}^{\circ}$ and $a_{l}=1_{i_l}$ for $l\neq k$, 
and otherwise $\Phi_{1^{n}}'(a_1, \ldots, a_n)=0$. 
\end{Corollary}
{\it Proof.}
This formula follows immediately from \eqref{eq:7.7}.
\hfill $\blacksquare$\\

\section{Higher order derivatives}

We can derive a formula for higher-order derivatives of 
$\Phi_{w,t}(a_1, \ldots, a_n)$ at $t=0$. We denote
{
\setlength{\belowdisplayshortskip}{\baselineskip}
\begin{equation}\label{eq:8.1}
\Phi^{(m)}_{w}(a_1,\ldots, a_n):=
\left.\left[\frac{d^{m}}{dt^{m}}\Phi_{w,t}(a_1,\ldots a_n)\right]\right|_{t=0}
\end{equation}
}%
for any $m\in {\mathbb N}$ and $w\in \mathpzc{M}_{n}$, where $a_1\in {\mathcal A}_{i_1}^{\circ}, 
\ldots, a_n\in {\mathcal A}_{i_n}^{\circ}$. 

The main feature of this formula is that 
at each local maximum associated with $\{k\}$ of $w$ we have a singleton block and thus the corresponding
moment $\varphi_{i_{k}}(a_k)$ vanishes. Since the derivatives 
$\varphi_{i_k}^{(n)}(a_k)$ do not have to vanish, we have to distribute
the derivative of order $k$ in such a way that to each singleton block we assign 
a derivative of order at least $1$ and to the remaining level return blocks of 
$\pi(w)$ we assign the remaining order.

The Leibniz rule used in this computation refers to the situation in which we have a monomial
of $p+r$ functions,
\[
F(g_1(x), \ldots, g_{r}(x))=\prod_{k=1}^{p}g_{k}(x)
\]
such that $g_1(x_0)=\cdots =g_{p}(x_0)$ for some $p\leq r$, 
but higher order derivatives of these $g_{k}$ at $x_0$
are not assumed to vanish, and no conditions on $g_{p+1}(x_0), \ldots, g_r(x_0)$ are given.
In that situation, non-zero contribution to the derivative of order $m$ of $F$ with respect to
$x$ is given by products in which a derivative of order at least one is assigned to
each $g_k$, and the remaining ones are distributed in an arbitrary way. Thus we obtain
\[
F^{(m)}(g_1(x_0), \ldots, g_{r}(x_0))=
\sum_{\stackrel{m_1+\cdots +m_r=m}{\scriptscriptstyle m_{k}\geq 1\;{\rm for}\;k\leq p}}
\frac{m!}{m_1! \cdots m_r!}
\prod_{s=1}^{p}
g_{s}^{(m_s)}(x_0)
\;
\prod_{s=p+1}^{r} 
g_s^{(m_s)}(x_0)
\]
for any $m\geq 1$.

\begin{Theorem}
Let $a_1\in {\mathcal A}_{i_1}^{\circ}, \ldots, a_n\in \mathcal{A}_{i_n}^{\circ}$, 
where $i_1\neq \cdots \neq i_n$ and let $w=j_1\cdots j_n\in \mathpzc{M}_{n}$.
Let $\pi(w)=\{V_1, \ldots, V_{r}\}$, where $V_1=\{k_1\}, \ldots, V_p=\{k_p\}$ are singleton blocks 
and $V_{p+1}, \ldots, V_r$ are other blocks, where $p\leq r$. 
\begin{enumerate}
\item 
If $p\leq m$, then it holds that
\begin{align*}
\;\;\;\;\;\;\Phi_{w}^{(m)}(a_1, \ldots, a_n)&=
\sum_{\stackrel{m_1+\cdots +m_r=m}{\scriptscriptstyle m_{k}\geq 1\;{\rm for}\;k\leq p}}
\frac{m!}{m_1! \cdots m_r!}
\prod_{s=1}^{p}
\varphi_{i_{k_s}}^{(m_s)}(a_{k_s})
\prod_{s=p+1}^{r} \beta_{|V_s|}^{(m_s)}((a_{k})_{k\in V_s})
\end{align*}
whenever $w$ is adapted to $(i_1, \ldots, i_n)$.
\item
If $w$ is not adapted to $(i_1, \ldots, i_n)$
or $p>m$, then $\Phi_{w}^{(m)}(a_1, \ldots, a_n)=0$.
\end{enumerate}
\end{Theorem}
{\it Proof.}
It is enough to apply the Leibniz formula given above to the computation of derivatives
of moments given by Theorem 6.1. The role of $g_{k}(x)$ for $1\leq k \leq p$ is played by 
$\varphi_{i_k,t}(a_{k})$ since $\varphi_{i_k,0}(a_k)=0$ 
and the role of $g_{k}(x)$ for $k>p$ is played by $\beta_{|V_s|}((a_{k})_{k\in V})$. 
Of course, we need to take into account only those words $w\in \mathpzc{M}_{n}$ 
which are adapted to $\ell=(i_1, \ldots, i_n)$ since only for these words the moment
functional $\Phi_{w,t}(a_1, \ldots , a_n)$ may not be zero.
\hfill $\blacksquare$

\begin{Corollary}
Let $\varphi=(\star_{i\in I}\varphi_{i})\circ \tau$ and let $\varphi^{(m)}$ 
denote its derivative of order $m$. Under the assumptions of Theorem 8.1, it holds that
\[
\varphi^{(m)}(a_1\cdots a_n)=
\sum_{\stackrel{w\in \mathpzc{M}_{n}(\ell)}{\scriptscriptstyle LM(w)\leq m}}
\Phi_{w}^{(m)}(a_1, \ldots, a_n)
\]
where $\mathpzc{M}_{n}(\ell)$ is the subset of $\mathpzc{M}_{n}$ consisting of 
words adapted to $\ell=(i_1, \ldots, i_n)$ and $LM(w)$ is the number of local maxima 
in $w$.
\end{Corollary}
{\it Proof.}
This assertion follows from \eqref{eq:4.8} and Theorem 8.1. 
\hfill $\blacksquare$

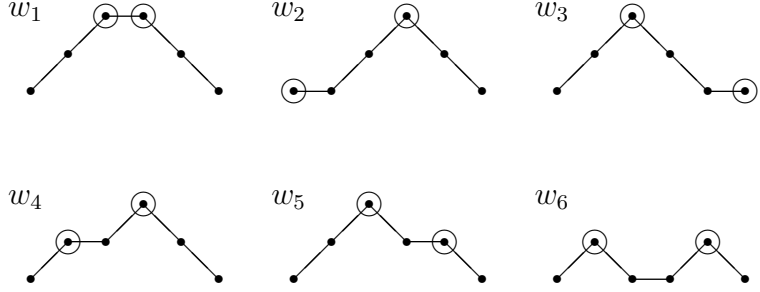
\begin{figure}
\unitlength=1mm
\special{em:linewidth 1pt}
\linethickness{0.5pt}
\begin{picture}(180.00,50.00)(-35.00,25.00)

%%%%%%%%%%%%%%%path %%%%%%%%%%%%%%%
\put(-8.00,65.00){$w_1$}
\put(27.00,65.00){$w_2$}
\put(62.00,65.00){$w_3$}
\put(-8.00,40.00){$w_4$}
\put(27.00,40.00){$w_5$}
\put(62.00,40.00){$w_6$}
%%%%%%%%%%%%%%%%% path 1%%%%%%%%%%%%%%%%%%%%%%%%

\put(-5.00,55.00){\line(1,1){5.00}}
\put(0.00,60.00){\line(1,1){5.00}}
\put(5.00,65.00){\line(1,0){5.00}}
\put(10.00,65.00){\line(1,-1){5.00}}
\put(15.00,60.00){\line(1,-1){5.00}}

\put(-5.00,55.00){\circle*{1.00}}
\put(0.00,60.00){\circle*{1.00}}
\put(5.00,65.00){\circle*{1.00}}
\put(5.00,65.00){\circle{3.00}}
\put(10.00,65.00){\circle*{1.00}}
\put(10.00,65.00){\circle{3.00}}
\put(15.00,60.00){\circle*{1.00}}
\put(20.00,55.00){\circle*{1.00}}

%%%%%%%%%%%  path 2 %%%%%%%%%%%%%%%%%%%%%%%

\put(30.00,55.00){\line(1,0){5.00}}
\put(35.00,55.00){\line(1,1){5.00}}
\put(40.00,60.00){\line(1,1){5.00}}
\put(45.00,65.00){\line(1,-1){5.00}}
\put(50.00,60.00){\line(1,-1){5.00}}

\put(30.00,55.00){\circle*{1.00}}
\put(30.00,55.00){\circle{3.00}}
\put(35.00,55.00){\circle*{1.00}}
\put(40.00,60.00){\circle*{1.00}}

\put(45.00,65.00){\circle*{1.00}}
\put(45.00,65.00){\circle{3.00}}
\put(50.00,60.00){\circle*{1.00}}
\put(55.00,55.00){\circle*{1.00}}

%%%%%%%%%%%  path 3 %%%%%%%%%%%%%%%%%%%%%%%

\put(65.00,55.00){\line(1,1){5.00}}
\put(70.00,60.00){\line(1,1){5.00}}
\put(75.00,65.00){\line(1,-1){5.00}}
\put(80.00,60.00){\line(1,-1){5.00}}
\put(85.00,55.00){\line(1,0){5.00}}

\put(65.00,55.00){\circle*{1.00}}
\put(70.00,60.00){\circle*{1.00}}
\put(75.00,65.00){\circle*{1.00}}
\put(75.00,65.00){\circle{3.00}}
\put(80.00,60.00){\circle*{1.00}}
\put(85.00,55.00){\circle*{1.00}}
\put(90.00,55.00){\circle*{1.00}}
\put(90.00,55.00){\circle{3.00}}

%%%%%%%%%%%%%%%%% path 4 %%%%%%%%%%%%%%%%%%%%%%%%
\put(-5.00,30.00){\line(1,1){5.00}}
\put(0.00,35.00){\line(1,0){5.00}}
\put(5.00,35.00){\line(1,1){5.00}}
\put(10.00,40.00){\line(1,-1){5.00}}
\put(15.00,35.00){\line(1,-1){5.00}}

\put(-5.00,30.00){\circle*{1.00}}
\put(0.00,35.00){\circle*{1.00}}
\put(0.00,35.00){\circle{3.00}}
\put(5.00,35.00){\circle*{1.00}}
\put(10.00,40.00){\circle*{1.00}}
\put(10.00,40.00){\circle{3.00}}
\put(15.00,35.00){\circle*{1.00}}
\put(20.00,30.00){\circle*{1.00}}

%%%%%%%%%%%  path 5 %%%%%%%%%%%%%%%%%%%%%%%

\put(30.00,30.00){\line(1,1){5.00}}
\put(35.00,35.00){\line(1,1){5.00}}
\put(40.00,40.00){\line(1,-1){5.00}}
\put(45.00,35.00){\line(1,0){5.00}}
\put(50.00,35.00){\line(1,-1){5.00}}

\put(30.00,30.00){\circle*{1.00}}
\put(35.00,35.00){\circle*{1.00}}
\put(40.00,40.00){\circle*{1.00}}
\put(40.00,40.00){\circle{3.00}}
\put(45.00,35.00){\circle*{1.00}}
\put(50.00,35.00){\circle*{1.00}}
\put(50.00,35.00){\circle{3.00}}
\put(55.00,30.00){\circle*{1.00}}

%%%%%%%%%%%  path 6 %%%%%%%%%%%%%%%%%%%%%%%

\put(65.00,30.00){\line(1,1){5.00}}
\put(70.00,35.00){\line(1,-1){5.00}}
\put(75.00,30.00){\line(1,0){5.00}}
\put(80.00,30.00){\line(1,1){5.00}}
\put(85.00,35.00){\line(1,-1){5.00}}

\put(65.00,30.00){\circle*{1.00}}
\put(70.00,35.00){\circle*{1.00}}
\put(75.00,30.00){\circle*{1.00}}
\put(70.00,35.00){\circle{3.00}}
\put(80.00,30.00){\circle*{1.00}}
\put(85.00,35.00){\circle*{1.00}}
\put(90.00,30.00){\circle*{1.00}}
\put(85.00,35.00){\circle{3.00}}

%%%% partition  %%%%%%%%%%%%%%%%%%%%%%

\end{picture}
\caption{Motzkin paths of length $5$ with two local maxima.}
\end{figure}

\begin{Example}
{\rm For $n=6$ we do not have a pyramid path, therefore in the computation of 
the second order derivatives $\Phi_{w}^{(2)}$ it is enough to take into account six
paths of length $5$ shown in Fig. 3 (paths of length $5$ give words of length $6$).
Each of these paths has two local maxima, therefore the formula of Theorem 8.1 reduces 
to one term. We obtain
{
\setlength{\belowdisplayshortskip}{\baselineskip}
\[
\Phi_{w}^{(2)}(a_1,a_2,a_3,a_4,a_5,a_6)=
\left\{
\begin{array}{lll}
2\varphi_{i_1}(a_1a_6)\varphi_{i_2}(a_2a_5)\varphi_{i_3}'(a_3)\varphi_{i_4}'(a_4)& {\rm if}&w=w_1\\[5pt]
2\varphi_{i_2}(a_2a_6)\varphi_{i_3}(a_3a_5)\varphi_{i_1}'(a_1)\varphi_{i_4}'(a_4)& {\rm if} &w=w_2\\[5pt]
2\varphi_{i_1}(a_1a_5)\varphi_{i_2}(a_2a_4)\varphi_{i_3}'(a_3)\varphi_{i_6}'(a_6)& {\rm if} &w=w_3\\[5pt]
2\varphi_{i_1}(a_1a_6)\varphi_{i_3}(a_3a_5)\varphi_{i_2}'(a_2)\varphi_{i_4}'(a_4)& {\rm if} &w=w_4\\[5pt]
2\varphi_{i_1}(a_1a_6)\varphi_{i_2}(a_2a_4)\varphi_{i_3}'(a_3)\varphi_{i_5}'(a_5)& {\rm if} &w=w_5\\[5pt]
2\varphi_{i_1}(a_1a_3)\varphi_{i_4}(a_4a_6)\varphi_{i_2}'(a_2)\varphi_{i_5}'(a_5)& {\rm if} &w=w_6\\[5pt]
\end{array}
\right.
\]
}%
where, for each $w_k$ separately, we assume that it is adapted to $(i_1, \ldots, i_6)$.
By Corollary 8.1, the sum of these expressions, with the appropriate Kronecker deltas, gives
{
\setlength{\belowdisplayshortskip}{\baselineskip}
\begin{align*}
\varphi^{(2)}(a_1a_2a_3a_4a_5a_6)&=\delta_{i_1,i_6}\delta_{i_2,i_5}\Phi_{w_1}^{(2)}+
\delta_{i_2,i_6}\delta_{i_3,i_5}\Phi_{w_2}^{(2)}
+\delta_{i_1,i_5}\delta_{i_2,i_4}\Phi_{w_3}^{(2)}\\
&\quad +\delta_{i_1,i_6}\delta_{i_3,i_5}\Phi_{w_4}^{(2)}
+\delta_{i_1,i_6}\delta_{i_2,i_4}\Phi_{w_5}^{(2)}+\delta_{i_1,i_3}\delta_{i_4,i_6}\Phi_{w_6}^{(2)}
\end{align*}
}%
where $\Phi_{w_k}^{(2)}:=\Phi_{w_k}^{(2)}(a_1, \ldots, a_6)$ 
and $\varphi=(\star_{i\in I}\varphi_{i})\circ \tau$
is the free product of functionals $\varphi_i$ composed with the unit identification map.
}
\end{Example}

\begin{Remark}
{\rm 
It is typical for the computations of the second derivative that non-zero contribution is obtained 
from three types of paths: flat pyramid with a plateau of length $1$ (of type $UHD$),
paths with one strong peak (of type $UD$) and one weak peak ($UH$ or $HD$), paths with two strong peaks.
It can be shown that the number of reduced Motzkin words of length $n$ which have two local maxima is equal to
\[
|\{w\in \mathpzc{M}_{n}: LM(w)=2\}|=
\left\{
\begin{array}{ll}
\binom{m}{2}& {\rm if}\;n=2m+1\\[5pt]
\binom{m+1}{2}& {\rm if}\;n=2m
\end{array}
\right.
\]
For instance, for $n=6$ we have 
$6$ paths shown in Fig. 3. The sequence obtained from the above formulas is given by: $0,1,0,3,1,6,3,10,6,15,10,21,15,\ldots $.
}
\end{Remark}

\section{Infinitesimal conditional freeness}

Using a similar approach to that presented for infinitesimal freeness, we are going to 
treat infinitesimal conditional freeness.

Let us begin with the definition of conditional freeness given in \cite{BLS}. Since it involves pairs of functionals, 
for instance $(\varphi, \psi)$, we will distinguish elements from ${\rm Ker}(\varphi)$ from
elements from ${\rm Ker}(\psi)$ by adopting special notations:
\vspace{3pt}
\begin{equation}\label{eq:9.1}
{\mathcal A}_{i}^{\circ}:={\mathcal A}_{i}\cap {\rm Ker}(\varphi)\;\;\;\;{\rm and}\;\;\;\;
{\mathcal A}_{i}^{\square}:={\mathcal A}_{i}\cap {\rm Ker}(\psi),
\end{equation}
it will be important to distinguish between them. A similar notation will be used for
localized pairs of functionals $(\varphi_{i}, \psi_{i})$, where $i\in I$.

\begin{Definition}
{\rm Let $(\varphi, \psi)$ be a pair of normalized linear functionals on a unital algebra ${\mathcal A}$.
The family $({\mathcal A}_{i})_{i\in I}$ of unital subalgebras of ${\mathcal A}$ is 
{\it conditionally free} (or, {\it c-free}) with respect to $(\varphi, \psi)$ if and only if
\begin{equation}\label{eq:9.2}
\varphi(a_1\cdots a_n)=\varphi(a_1)\cdots \varphi(a_n)
\end{equation}
whenever $a_i\in {\mathcal A}_{i}^{\square}:={\mathcal A}_{i}\cap {\rm Ker}(\psi)$ for all $i\in [n]$
and $i_1\neq \cdots \neq i_n$.}
\end{Definition}

It is convenient to use an equivalent definition which is more similar to that of freeness given by \eqref{eq:2.1}.
We state this as a proposition and provide its elementary proof.\\

\begin{Proposition}
The family $({\mathcal A}_{i})_{i\in I}$ of unital subalgebras of a unital algebra 
${\mathcal A}$ is conditionally free with respect to $(\varphi, \psi)$ if and only if 
\begin{equation}\label{eq:9.3}
\varphi(a_1\cdots a_n)=0
\end{equation}
whenever $a_1\in {\mathcal A}_{1}^{\circ}:={\mathcal A}_{i}\cap {\rm Ker}(\varphi)$
and $a_k \in {\mathcal A}_{i_k}^{\square}={\mathcal A}_{i_k}\cap {\rm Ker}(\psi)$ 
for all $i>1$, where $i_1\neq \cdots \neq i_n$. 
\end{Proposition}
{\it Proof.}
The family $({\mathcal A}_{i})_{i\in I}$ is c-free with respect to $(\varphi, \psi)$ if and only if 
\[
\varphi(a_1^{\square}\cdots a_n^{\square})=(\varphi(a_1)-\psi(a_1))\cdots (\varphi(a_n)-\psi(a_n))
\]
whenever $a_{k}^{\square}=a_{k}-\psi(a_k)$ for all $i\in [n]$ where $a_k\in {\mathcal A}_{i_k}$ 
for all $k$ and $i_1\neq \cdots \neq i_n$ and otherwise are arbitrary. 
This holds if and only if
\[
\varphi(a_1a_{2}^{\square}\cdots a_{n}^{\square})
=\varphi(a_1)(\varphi(a_2)-\psi(a_2))\cdots (\varphi(a_n)-\psi(a_n))
\]
for any $a_k\in {\mathcal A}_{i_k}$, $k\in [n]$, where $i_1\neq \cdots \neq i_n$, 
which is equivalent to
\[
\varphi(a_1^{\circ}a_{2}^{\square}\cdots a_{n}^{\square})=0
\]
where $a_{1}^{\circ}=a_1- \varphi(a_1) \in {\rm Ker}(\varphi)$. 
Thus we obtain \eqref{eq:9.3}, with the role of 
$a_1, a_2, \ldots, a_n$ played by $a_1^{\circ}, a_2^{\square}, \ldots, a_n^{\square}$.
This completes the proof. 
\hfill $\blacksquare$\\

\begin{Definition}
{\rm 
Let $(\varphi_{i}, \psi_{i})_{i\in I}$ be a family of pairs of normalized linear functionals 
on unital algebras ${\mathcal A}_{i\in I}$, respectively.
As in the case of infinitesimal freeness, we introduce a continuous deformation 
of these families
\begin{align*}
\varphi_{i,t}&=\varphi_{i}+\frac{\varphi_{i}^{(1)}}{1!}t+\frac{\varphi_{i}^{(2)}}{2!}t^2+\cdots +
\frac{{\varphi}_{i}^{(m)}}{m!}t^{m}+o(t^{m})\\
\psi_{i,t}&=\psi_{i}+\frac{\psi_{i}^{(1)}}{1!}t+\frac{\psi_{i}^{(2)}}{2!}t^2+\cdots +
\frac{{\psi}_{i}^{(m)}}{m!}t^{m}+o(t^{m})
\end{align*}
where $(\varphi_{i}^{(k)})_{k\geq 1}$ and $(\psi_{i}^{(k)})_{k\geq 1}$ are, for each $i\in I$, additional families 
of linear functionals on ${\mathcal A}_{i}$, such that $\varphi_{i}^{(k)}(1_i)=\psi_{i}^{(k)}(1_i)=0$
for all $k\geq 1$, where $m\in {\mathbb N}$ and $t\in (-\epsilon, \epsilon)$.
}
\end{Definition}

Note that the functionals $\varphi_{i,t}$ and $\psi_{i,t}$ are
normalized for any $i\in I$ and any $t$. 
As in the case of infinitesimal freeness, we will write
as $\varphi_{i}^{(1)}=\varphi_{i}'$ and $\psi_{i}^{(1)}=\psi_i'$ since these functionals
can be treated as first-order derivatives of $\phi_{i,t}$ and $\psi_{i,t}$ at $t=0$, respectively.
For these continuous deformations, we define their Boolean extensions  
$\widetilde{\varphi}_{i,t}$ and $\widetilde{\psi}_{i,t}$, respectively, by the formulas
\begin{align*}
\widetilde{\varphi}_{i,t}(p_{i}^{\alpha}a_1p_i\cdots p_ia_np_i^{\beta})&=
\varphi_{i,t}(a_1)\cdots \varphi_{i,t}(a_n)\\[3pt]
\widetilde{\psi}_{i,t}(p_{i}^{\alpha}a_1p_i\cdots p_ia_np_i^{\beta})&=
\psi_{i,t}(a_1)\cdots \psi_{i,t}(a_n)
\end{align*}
for any $a_1, \ldots, a_n\in {\mathcal A}_{i}$ and $\alpha, \beta\in \{0,1\}$, where $(p_i)_{i\in I}$ 
is a family of idempotents.
Equivalently, we can say that we use a simple version of infinitesimal freeness of Mingo and Tseng \cite{MT}
localized to each ${\mathcal A}_{i}$ separately to construct the 
extended functionals $\widetilde{\varphi}_{i,t}$ on $\widetilde{\mathcal A}_{i}$, respectively.

Using these families, we can define infinitesimally deformed tensor product functionals which are 
the c-free counterparts of $\Phi_{\otimes, t}$ given by \eqref{eq:4.10}. We need two functionals, of which the first one plays the 
role of a truly new functional, from which we obtain the first functional in the c-free product framework,
whereas the second one is of the same form as $\Phi_{\otimes, t}$, but is built only from the functionals $\psi_{i,t}$.
Therefore, let 
\vspace{3pt}
{
\setlength{\belowdisplayshortskip}{\baselineskip}
\begin{equation*}
\widetilde{\Phi}_{\otimes, t}=\bigotimes_{i\in I}(\widetilde{\varphi}_{i,t}\otimes \widetilde{\psi}_{i,t}^{\otimes \infty})\;\;\;{\rm and}\;\;\;
\widetilde{\Psi}_{\otimes, t}=\bigotimes_{i\in I}\widetilde{\psi}_{i,t}^{\otimes \infty}
\end{equation*}
}%
for any $t$. In other words, to construct the first functional, we apply a continuous deformation and an infinitesimal idempotent to 
each label and color as before, except that we use a different family for color $1$ and a different one for the 
remaining colors. 

From the infinitesimal perspective, the structure of our construction is based on two types of infinitesimal operations:
a continuous one (implemented by $t$) and operatorial ones (implemented by the
infinitesimal idempotents $p_i$):
{
\setlength{\belowdisplayshortskip}{\baselineskip}
\[
\varphi_{i}\xrightarrow {\text{deform}}
\varphi_{i,t}\xrightarrow {\text{extend}}
\widetilde{\varphi}_{i,t}
\]
\[
\psi_{i}\xrightarrow{\text{deform}} \psi_{i,t}\xrightarrow{\text{extend}}\widetilde{\psi}_{i,t}
\]
}%
\noindent
for any $i\in I$ and any $t$. These two types of local infinitesimal 
operations give the globally deformed functionals $\widetilde{\Phi}_{\otimes,t}$ and $\widetilde{\Psi}_{\otimes, t}$,
from which we extract their restrictions to the replica space ${\mathcal A}$, namely 
$\widetilde{\Phi}_{t}$ and $\widetilde{\Psi}_{t}$.
Then we define multilinear functionals $\widetilde{\Phi}_{w,t}$ and $\widetilde{\Psi}_{w,t}$ 
as in Definition 4.5 and compute their derivatives at $t=0$. Since tensor products admit canonical
derivations, the computation of these derivatives is again natural.

\begin{Definition}
{\rm For a family of pairs $(\varphi_{i}, \psi_{i})_{i\in I}$ of normalized functionals 
on ${\mathcal A}_{i}$, respectively, define multilinear functionals by the formulas
\vspace{3pt}
{
\setlength{\belowdisplayshortskip}{\baselineskip}
\begin{align}\label{eq:9.4}
\widetilde{\Phi}_{w,t}(a_1, \ldots, a_n)&=\;\;\widetilde{\Phi}_{t}(a_1(j_1)\cdots a_n(j_n))\\
\label{eq:9.5}
\widetilde{\Phi}_{w}'\;(a_1, \ldots, a_n)\;&=\left.\left[\frac{d\widetilde{\Phi}_{w,t}(a_1, \ldots, a_n)}{dt}\right]\right|_{t=0}
\end{align}
}%
where $w=j_1\cdots j_n\in \mathpzc{M}_{n}$, $n\in {\mathbb N}$ and 
$a_1\in {\mathcal A}_{i_1}, \ldots , a_n\in {\mathcal A}_{i_n}$, $i_1\neq \cdots \neq i_n$.
Multilinear functionals $\widetilde{\Psi}_{w,t}$ and $\widetilde{\Psi}_{w}'$ are defined in an analogous fashion. 
}
\end{Definition}

As in the case of $\Phi_{t}$ and $\Phi_{w}$, 
we slightly abuse notation, since we consider linear functionals $\widetilde{\Phi}_t$ 
indexed by a continuous parameter $t$ as well as multilinear functionals $\widetilde{\Phi}_w$ indexed by
Motzkin paths $w$ (we do the same for $\widetilde{\Psi}_t$ and $\widetilde{\Psi}_{w}$). 
This should not cause confusion because $t$ and $w$ are entirely different objects. 

The main motivation to study such deformations of $\widetilde{\Phi}$ is the fact that the 
functional $\widetilde{\Phi}$ is the c-free analog of $\Phi$ and reproduces moments of c-free random variables
as shown in \cite{L1}. We observed in \cite{L4} that the decomposition of the free product
of functionals given by Theorem 4.1 can be generalized to the conditionally free product 
of functionals $*_{i\in I}(\varphi_{i}, \psi_{i})$. Then we can differentiate this 
product in the same way as we did in Theorem 7.1. Below we give an explicit formula 
for the derivative $\varphi'$, the first functional from the pair $(\varphi, \psi)$. 
The derivative $\psi'$ is computed in the same way as in Theorem 7.1 since 
$\psi$ is the free product of the $\psi_{i,t}$ and that is why it will be omitted.

\begin{Proposition}
Let $({\mathcal A}_{i},\varphi_i)_{i\in I}$ be noncommutative probability spaces. 
Let $\varphi'=d\varphi_{t}/dt|_{t=0}$ and $\psi'=d\psi_{t}/dt|_{t=0}$, where 
$(\varphi_{t}, \psi_{t})=(\star_{i\in I} (\varphi_{i,t}, \psi_{i,t}))\circ \tau$.
Let $a_1\in \mathcal{A}_{i_1}, \ldots, a_n\in \mathcal{A}_{i_n}$, where $i_1\neq \cdots \neq i_n$, $n\in {\mathbb N}$.
Then the following decomposition holds:
\begin{equation}\label{eq:9.6}
\varphi'(a_1\cdots a_n)=
\sum_{w\in \mathpzc{M}_{n}}
\widetilde{\Phi}_{w}'(a_{1}, \ldots, a_n)
\end{equation}
The decomposition of $\psi'$ is given by Theorem 7.1, with the $\psi_{i}$ replacing the $\varphi_{i}$.
\end{Proposition}
{\it Proof.}
The decomposition of $\varphi$ in terms of $\widetilde{\Phi}_{w}$ and thus 
the decomposition of $\varphi_{t}$ in terms of $\widetilde{\Phi}_{w,t}$ is completely analogous to that of 
Theorem 4.1 and Corollary 4.1, respectively. This follows from the tensor product construction given in \cite{L1} which holds
for conditional freeness. In particular, the orthogonal replicas are defined in the same way and the only difference 
is that instead of $\Phi_{\otimes,t}$, we use  $\widetilde{\Phi}_{\otimes, t}$. This proves \eqref{eq:9.6}.
The decomposition of $\psi'$ is analogous to that in Theorem 7.1. This completes the proof.
\hfill $\blacksquare$\\

We need to compute $\widetilde{\Phi}_{w}'$ for $w\in \mathpzc{M}$. 
In our notations, we need to distinguish moments and cumulants 
associated with $\varphi_{i}$ from those associated with $\psi_{i}$.
Therefore, let 
\begin{align*}
\beta_{|V|}^{\,\chi}((a_k)_{k\in V})&=\beta_{|V|}^{\,\chi}(a_{k_1}, \ldots , a_{k_p})
\end{align*}
be the Boolean cumulants associated with a functional  $\chi$,
where $V=\{k_1<\cdots <k_p\}$. 
We will then adopt the following joint notation for moment functionals:
\vspace{3pt}
\[
\eta_{i}=\left\{
\begin{array}{ll}
\varphi_{i} & {\rm if}\;\; j_i=1\\
\psi_{i} & {\rm if}\;\; j_i>1
\end{array}
\right.
\]
and Boolean cumulant functionals:
\vspace{3pt}
\[
\gamma_{|V|}=\left\{
\begin{array}{ll}
\beta_{|V|}^{\,\varphi} & {\rm if}\;\; j_{k_1}=\cdots=j_{k_p}=1\\
\beta_{|V|}^{\,\psi}   & {\rm if}\;\; j_{k_1}=\cdots=j_{k_p}>1 
\end{array}
\right.
\]\\
In terms of the geometry of Motzkin paths, 
this means that moments and cumulants of variables $a_k$ corresponding to the 
lowest level are associated with $\varphi_{i_k}$ and those corresponding to 
the remaining levels are associated with $\psi_{i_k}$.
\vspace{3pt}

\begin{Proposition}
Let $a_1\in {\mathcal A}_{i_1}, \ldots, a_n\in \mathcal{A}_{i_n}$, where $i_1\neq \cdots \neq i_n$
and let $w\in \mathpzc{M}_{n}$. Then 
\begin{equation}\label{eq:9.7}
\widetilde{\Phi}_{w}'(a_1, \ldots, a_n)=\sum_{V\in \pi(w)}\gamma_{|V|}'((a_{i})_{i\in V})
\prod_{U\neq V}\gamma_{|U|}((a_{i})_{i\in U})
\end{equation}
whenever $w$ is adapted to $(i_1, \ldots, i_n)$ and otherwise it vanishes.
\end{Proposition}
{\it Proof.}
The only difference between the former deformed Motzkin functionals $\Phi_{w,t}$ and their c-free
counterparts $\widetilde{\Phi}_{w,t}$ is that tensor positions at level $1$ are occupied by functionals 
$\widetilde{\psi}_{i,t}$ instead of $\widetilde{\varphi}_{i,t}$ and that is taken care of 
by the new two-functional notation. We have 
\[
\widetilde{\Phi}_{w,t}(a_1, \ldots, a_n)=
\prod_{V\in \pi(w)}\gamma_{|V|,t}((a_{i})_{i\in V})
\]
whenever $w$ is adapted to $(i_1, \ldots, i_n)$, where $\pi(w)$ is the level return partition associated with $w$,
as in \eqref{eq:6.1}, where $\gamma_{|V|,t}$ stands for the $t$-deformation of 
$\gamma_{|V|}$ given by
\[
\gamma_{|V|,t}=\left\{
\begin{array}{ll}
\beta_{|V|}^{\,\varphi_{t}} & {\rm if}\;\; j_{k_1}=\cdots=j_{k_p}=1\\
\beta_{|V|}^{\,\psi_{t}}   & {\rm if}\;\; j_{k_1}=\cdots=j_{k_p}>1 
\end{array}
\right.
\]\\
If $w$ is not adapted to $(i_1, \ldots, i_n)$, then $\widetilde{\Phi}_{w,t}(a_1, \ldots, a_n)=0$.  
Then it suffices to use the Leibniz rule, as in \eqref{eq:7.3}, which completes our proof.
\hfill $\blacksquare$\\

So far the formulas look very similar to those in the free case, but the difference is hidden 
in the notation. In fact, it turns out two types of Motzkin paths play a special role in the infinitesimal setting: flat paths and concatenations of a pyramid path and a flat path (including the case when this flat path is empty if $n$ is odd).

\begin{Theorem}
Let $a_1\in {\mathcal A}_{i_1}^{\circ}, a_2\in {\mathcal A}_{i_2}^{\square}\ldots, a_n\in \mathcal{A}_{i_n}^{\square}$, 
where $i_1\neq \cdots \neq i_n$ and let $w=j_1\cdots j_n\in \mathpzc{M}_{n}$.
\begin{enumerate}
\item If $w=1^n$ for $n\in {\mathbb N}$, then 
\begin{equation}\label{eq:9.8}
\widetilde{\Phi}_{w}'(a_1, \ldots, a_n)=\varphi_{i_1}'(a_1)\varphi_{i_2}(a_2)\cdots \varphi_{i_n}(a_n).
\end{equation}
\item If $w=w_1w_2$, where $w_1$ is a pyramid 
adapted to $(i_1, \ldots, i_{2m-1})$ and $w_2=1^{n-2m+1}$, where $2m-1\leq n$, then
\begin{equation}\label{eq:9.9}
\begin{split}
\widetilde{\Phi}_{w}'(a_1, \ldots, a_n)&=\varphi_{i_1}(a_1a_{2m-1})\varphi_{i_{2m}}(a_{2m})\cdots \varphi_{i_n}(a_{n})\\
&\quad \times\;\psi_{i_2}(a_2a_{2m-2})\cdots \psi_{i_{m-1}}(a_{m-1}a_{m+1})\psi_{i_m}'(a_m).
\end{split}
\end{equation}
\item
In the remaining cases, $\widetilde{\Phi}_{w}'(a_1, \ldots, a_n)=0$.
\end{enumerate}
\end{Theorem}
{\it Proof.}
The case $n=1$ is obvious, so let us assume that $n>1$. 
As in the case of $\Phi_w'$, we need to identify local maxima since they produce singletons. However, in contrast 
to singletons in the infinitesimally free case, not all singletons have to be taken into account since we assume that 
\[
\psi_{i_2}(a_2)=\ldots =\psi_{i_n}(a_n)=0
\] 
and if $a_n$ is located at level $1$, the 
computation of moments $\widetilde{\Phi}_w'(a_1, \ldots a_n)$ gives $\varphi_{i_k}(a_k)$ for $k>1$ which does not have to vanish. Therefore, the only local maxima to which we need to assign a derivative are those indices $k>1$ to which 
corresponds letter $j_k>1$ (then we need to assign to it $\psi_{i_k}'(a_k)$) or 
the index $k=1$ to which corresponds letter $j_1=1$ (then we need to assign to it $\varphi_{i_1}'(a_1)$). 
Since in the computation of the first derivative w.r.t. $t$ only one `local' derivative is available, 
the only words $w$ which may give a non-zero contributions are those which correpond
to: 1) a concatenation of a pyramid that begins at $1$ and ends at $2m-1\leq n$ for some $m$ and a flat path 
from $2m-1$ to $n$, or a flat path from $1$ to $n$. The path of the first type has one local maximum at $(m,m)$ 
which contributes $\psi_{i_m}'(a_m)$, with the remaining factors given by Boolean cumulants associated with 
the level return blocks obtained from the excursion from $1$ to $2m-1$ and factors of type $\varphi_{i_k}(a_k)$
for $k=2m, \ldots, n$. Therefore, the moment associated with such a concatenation is 
of the form
\begin{align*}
\widetilde{\Phi}_{w}'(a_1, \ldots, a_n)&=\varphi_{i_1}(a_1a_{2m-1})\varphi_{i_{2m}}(a_{2m})\cdots \varphi_{i_n}(a_{n})\\
&\quad  \psi_{i_2}(a_2a_{2m-2})\cdots \psi_{i_{m-1}}(a_{m-1}a_{m+1})\psi_{i_m}'(a_m).
\end{align*}
In turn, if $w$ is a constant word $1^{n}$, then we must assign to the first index the 'local' 
derivative $\varphi_{i_1}'(a_1)$ and the remaining factors are $\varphi_{i_2}(a_2), \ldots, \varphi_{i_n}(a_n)$, which gives
\[
\widetilde{\Phi}_{w}'(a_1, \ldots, a_n)=\varphi_{i_1}'(a_1)\varphi_{i_2}(a_2)\cdots \varphi_{i_n}(a_n).
\]
This completes the proof.
\hfill $\blacksquare$\\

\begin{figure}
\unitlength=1mm
\special{em:linewidth 1pt}
\linethickness{0.5pt}
\begin{picture}(180.00,35.00)(-35.00,45.00)

%%%%%%%%%%%%%%%%% path 1%%%%%%%%%%%%%%%%%%%%%%%%

\put(-5.00,55.00){\line(1,0){5.00}}
\put(0.00,55.00){\line(1,0){5.00}}
\put(5.00,55.00){\line(1,0){5.00}}
\put(10.00,55.00){\line(1,0){5.00}}

\put(-5.00,55.00){\circle*{1.00}}
\put(0.00,55.00){\circle*{1.00}}
\put(5.00,55.00){\circle*{1.00}}
\put(-5.00,55.00){\circle{3.00}}
\put(10.00,55.00){\circle*{1.00}}
\put(15.00,55.00){\circle*{1.00}}

%%%%%%%%%%%  path 2 %%%%%%%%%%%%%%%%%%%%%%%

\put(30.00,55.00){\line(1,1){5.00}}
\put(35.00,60.00){\line(1,-1){5.00}}
\put(40.00,55.00){\line(1,0){5.00}}
\put(45.00,55.00){\line(1,0){5.00}}

\put(30.00,55.00){\circle*{1.00}}
\put(35.00,60.00){\circle{3.00}}
\put(35.00,60.00){\circle*{1.00}}
\put(40.00,55.00){\circle*{1.00}}
\put(45.00,55.00){\circle*{1.00}}
\put(50.00,55.00){\circle*{1.00}}

%%%%%%%%%%%  path 3 %%%%%%%%%%%%%%%%%%%%%%%

\put(65.00,55.00){\line(1,1){5.00}}
\put(70.00,60.00){\line(1,1){5.00}}
\put(75.00,65.00){\line(1,-1){5.00}}
\put(80.00,60.00){\line(1,-1){5.00}}

\put(65.00,55.00){\circle*{1.00}}
\put(70.00,60.00){\circle*{1.00}}
\put(75.00,65.00){\circle*{1.00}}
\put(75.00,65.00){\circle{3.00}}
\put(80.00,60.00){\circle*{1.00}}
\put(85.00,55.00){\circle*{1.00}}

\end{picture}
\caption{Motzkin paths of length $4$ contributing to infinitesimal c-free moments. Local maxima of special type
are marked with hollow circles.}
\end{figure}
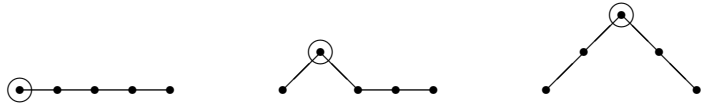

\begin{Example}
{\rm Let us compute the moments $\widetilde{\Phi}_{w}':=\widetilde{\Phi}_w'(a_1,a_2,a_3,a_4,a_5)$ corresponding to the Motzkin paths shown in Fig. 4. 
We assume that $a_1$ is $\varphi_{i_1}$-centered and each $a_k$ is $\psi_{i_k}$-centered for $k=2,3,4,5$.
Notice that these are all paths of length $4$ (words are of length $5$) which have non-zero derivatives provided 
each $w$ is adapted to $(i_1,\ldots, i_5)$. We obtain
{
\setlength{\belowdisplayshortskip}{\baselineskip}
\[
\widetilde{\Phi}_{w}'=
\left\{
\begin{array}{lll}
\varphi_{i_1}'(a_1)\varphi_{i_2}(a_2)\varphi_{i_3}(a_3)\varphi_{i_4}(a_4)\varphi_{i_5}(a_5)& {\rm if}& w=1^5\\[5pt]
\psi_{i_2}'(a_2)\varphi_{i_1}(a_1a_3)\varphi_{i_4}(a_4)\varphi_{i_5}(a_5)& {\rm if}& w= 121^3\\[5pt]
\psi_{i_3}'(a_3)\varphi_{i_1}(a_1a_5)\varphi_{i_2}(a_2a_4)& {\rm if}& w=12321
\end{array}
\right.
\]
}%
For instance, if we take $w=1^5$, then the derivative must be applied to 
$\varphi_{i_1}(a_1)$ in order to get a non-zero contribution. If we applied it
to $\varphi_{k}(a_k)$ for $k>1$, then we would have to keep $\varphi_{i_1}(a_1)$, but that one vanishes. 
Similarly, if we take $w=121^3$, then the derivative has to be applied to $\psi_{i_3}(a_3)$ since otherwise 
we would have to keep $\psi_{i_3}(a_3)$ which vanishes. Note also that in the computations
Boolean cumulants reduce to moments. For instance, for $w=12321$ we have 
$\gamma_{2}(a_1,a_5)=\varphi_{i_1}(a_1a_4)$ since $a_1$ is $\varphi_{i_1}$-centered and 
$\gamma_{2}(a_2,a_4)=\psi(a_2a_4)$ since $a_2, a_4$ are $\psi_{i_2}$-centered.
The sum of the above terms gives the following contribution to the derivative of the conditionally free product 
$(\varphi, \psi)=(\star_{i\in I}(\varphi_i, \psi_i))\circ \tau$:
{
\setlength{\belowdisplayshortskip}{\baselineskip}
\[
\varphi'(a_1,a_2,a_3,a_4,a_5)=\widetilde{\Phi}_{1^{5}}'+
\delta_{i_1,i_3}\widetilde{\Phi}_{121^{3}}'+
\delta_{i_1,i_5}\delta_{i_2,i_4}\widetilde{\Phi}_{12321}'
\]
}%
with 
{
\setlength{\belowdisplayshortskip}{\baselineskip}
\[
\psi'(a_1,a_2,a_3,a_4,a_5)=\delta_{i_1,i_5}\delta_{i_2,i_4}\widetilde{\Psi}_{12321}'
\]
}%
whenever $i_1\neq \cdots \neq i_5$. The Kronecker deltas reflect the adaptedness conditions within 
level return blocks. They are of the same type as in the case of a pyramid path except that 
they apply to the partial pyramid paths present in the given expression.
}
\end{Example}

We are ready to give a Leibniz-type definition of infinitesimal conditional freeness. We will show below that
this definition is motivated by Theorems 9.2 and 9.3.

\begin{Definition}
{\rm Let $({\mathcal A}, \varphi, \varphi', \psi,\psi')$ be a unital algebra equipped with 
two pairs of linear functionals, such that $\varphi(1)=\psi(1)=1$ and $\varphi'(1)=\psi'(1)=0$.
A family of unital subalgebras $({\mathcal A}_{i})_{i\in I}$ is 
{\it infinitesimally conditionally free} 
with respect to the quadruple $(\varphi, \varphi', \psi, \psi')$ if and only if
\begin{enumerate}
\item for any $a_1\in {\mathcal A}_{i_1}^{\circ}$, $a_2\in {\mathcal A}_{i_2}^{\square}, \ldots, a_n\in {\mathcal A}_{i_n}^{\square}$
it holds that 
\begin{equation}\label{eq:9.10}
\varphi'(a_1\cdots a_n)=\varphi'(a_1)\varphi(a_2 \cdots a_n)+
\sum_{m=2}^{n}\psi'(a_{m})\varphi(a_1\cdots a_{m-1}a_{m+1}\cdots a_n)
\end{equation}
whenever $i_1\neq \cdots \neq i_{n}$,
\item it is conditionally free with respect to $(\varphi, \psi)$,
\item it is infinitesimally free with respect to $(\psi, \psi')$.
\end{enumerate}
}
\end{Definition}
\vspace{3pt}

\begin{Theorem}
Under the assumptions of Definition 9.4, the family of unital subalgebras $({\mathcal A})_{i\in I}$
is infinitesimally conditionally free with respect to $(\varphi, \varphi', \psi, \psi')$ if and only if
\begin{enumerate}
\item
for any $a_1\in {\mathcal A}_{i_1}^{\circ}$, $a_2\in {\mathcal A}_{i_2}^{\square}, \ldots, a_n\in {\mathcal A}_{i_n}^{\square}$,
it holds that
\begin{equation}\label{eq:9.11} 
\begin{split}
\;\;\;\;\;\;\;\;\;\varphi'(a_1\cdots a_n)&=\varphi'(a_1)\varphi(a_2)\cdots \varphi(a_n) + \sum_{m=2}^{\lfloor n/2\rfloor}
\delta_{i_1,i_{2m-1}}\cdots \delta_{i_{m-1},i_{m+1}}\psi'(a_m)\\
&\!\!\quad \times
\varphi(a_1a_{2m-1})\psi(a_2a_{2m-2})\cdots \psi(a_{m-1}a_{m+1})\varphi(a_{2m})\cdots \varphi(a_{n}) 
\end{split}
\end{equation}
whenever $i_1\neq \cdots \neq i_n$,
\item
it is conditionally free with respect to $(\varphi, \psi)$,
\item
it is infinitesimally free with respect to $(\psi, \psi')$.
\end{enumerate}
\end{Theorem}
{\it Proof.}
It suffices to show the equivalence of (1) with Definition 9.4 (1). This is equivalent to showing that if
$a_1\in {\mathcal A}_{i_1}^{\circ}$, $a_2\in {\mathcal A}_{i_2}^{\square}, \ldots, a_n\in {\mathcal A}_{i_n}^{\square}$, where $i_1\neq \cdots \neq i_n$, then 
\begin{align*}
\varphi(a_1\cdots a_{m-1}a_{m+1}\cdots a_n)&=
\varphi(a_1a_{2m-1})\psi(a_2a_{2m-2})\cdots \psi(a_{m-1}a_{m+1})\varphi(a_{2m})\cdots \varphi(a_{n}) 
\end{align*}
if $2\leq m \leq \lfloor n/2\rfloor $ and $i_{1}=i_{2m-1}, \; \ldots, i_{m-1}=i_{m+1}$, and in the remaining cases 
this moment vanishes (in particular, if $m>\lfloor n/2\rfloor $). 
Suppose first that $2\leq m \leq \lfloor n/2 \rfloor $ and 
$i_{1}=i_{2m-1}, \; \ldots, i_{m-1}=i_{m+1}$. Then we can write 
\[
a_{m-1}a_{m+1}=(a_{m-1}a_{m+1})^{\square}+\psi_{i_{m-1}}(a_{m-1}a_{m+1})1
\]
and thus
\begin{align*}
\varphi(a_1\cdots a_{m-1}a_{m+1}\cdots a_{m})&=
\varphi(a_1\cdots (a_{m-1}a_{m+1})^{\square}\cdots a_{m})\\
&\quad +\psi(a_{m-1}a_{m+1})\varphi(a_1\cdots a_{m-2}a_{m+2}\cdots a_n)\\
&=\psi(a_{m-1}a_{m+1})\varphi(a_1\cdots a_{m-2}a_{m+2}\cdots a_n)
\end{align*}
by conditional freeness w.r.t. $\varphi$. It can be seen that an induction argument 
can be used and, in effect, it will reduce the expression to 
\[
\psi(a_2a_{2m-2})\cdots \psi(a_{m-1}a_{m+1})\varphi(a_1a_{2m-1}a_{2m}\cdots a_{n})
\]
Now, we use the decomposition
\[
a_1a_{2m-1}=(a_1a_{2m-1})^{\circ}+ \varphi(a_1a_{2m-1})1
\]
which gives 
\begin{align*}
\varphi(a_1a_{2m-1}a_{2m}\cdots a_{n})&=\varphi(a_1a_{2m-1})\varphi(a_{2m}\cdots a_{n})\\
&=\varphi(a_1a_{2m-1})\varphi(a_{2m})\cdots \varphi(a_{n})
\end{align*}
where we used the traditional definition of c-freeness w.r.t. $(\varphi,\psi)$.
This proves the required formula. In turn, if $i_{m-1}\neq i_{m+1}$, then all neighboring variables belong to different algebras
and therefore 
\[
\varphi(a_1\cdots a_{m-1}a_{m+1}\cdots a_{m})=0
\]
Finally, suppose that $m>\lfloor n/2 \rfloor $. Then we can reduce this mixed moment in a similar manner by succesive 
operations, but we will move with this reduction to the last variable $a_n$ and that is why we will obtain
{
\setlength{\belowdisplayshortskip}{\baselineskip}
\[
\varphi(a_1\cdots a_{m-1}a_{m+1}\cdots a_n)=\psi(a_{m-1}a_{m+1})\cdots \psi(a_{2m-n}a_{n})\varphi(a_1\cdots a_{2m-n-1})=0
\]
}%
\noindent
but this expression vanishes by Definition 9.1. This completes the proof.
\hfill $\blacksquare$\\

\begin{Corollary}
The formula for moments given by Theorem 9.2 can be written in the path decomposition form
\[
\varphi'(a_1\cdots a_n)=\widetilde{\Phi}_{1^n}'(a_1, \ldots, a_n)+\sum_{\stackrel{w=w_1w_2\in \mathpzc{M}_{n}}{\scriptscriptstyle w_1\;{\rm pyramid}, \;w_2\;{\rm flat}}}
\delta_{w_1}\widetilde{\Phi}_{w_1}'(a_{w_1})\widetilde{\Phi}_{w_2}(a_{w_2})
\]
where the sum ranges over all factorizations of $w$ such that $w_1$ is a non-empty pyramid path and $w_2$ is a flat path
and $\delta_{w_1}=1$ if $w_1$ is adapted to labels $i_1, \ldots, i_{|w_1|}$ and otherwise $\delta_{w_1}=0$.
\end{Corollary}
{\it Proof.}
This assertion follows from the formula \eqref{eq:9.11}. Let us add that in the summation 
we allow $w_2$ to be an empty path.
\hfill $\blacksquare$\\

\begin{Example}
{\rm Let us compute examples of infinitesimal moments of c-free random variables with the use 
of Theorem 9.2. In both examples, we assume that the first variable $a_1$ is $\varphi$-centered and the remaining 
ones are $\psi$-centered. 
We obtain
{
\setlength{\belowdisplayshortskip}{\baselineskip}
\begin{align*}
\varphi'(a_1a_2a_3)&= \varphi'(a_1)\varphi(a_2)\varphi(a_3)+\delta_{i_1,i_3}\psi'(a_2)\varphi(a_1a_3)\\
\varphi'(a_1a_2a_3a_4)&= \varphi'(a_1)\varphi(a_2)\varphi(a_3)\varphi(a_4)+
\delta_{i_1,i_3}\varphi(a_1a_3)\psi'(a_2)\varphi(a_4)\\
\varphi'(a_1a_2a_3a_4a_5)&= \varphi'(a_1)\varphi(a_2)\varphi(a_3)\varphi(a_4)\varphi(a_5)+
\delta_{i_1,i_3}\varphi(a_1a_3)\psi'(a_2)\varphi(a_4)\varphi(a_5)\\
&\quad +
\delta_{i_1,i_5}\delta_{i_2,i_4}\varphi(a_1a_5)\psi(a_2a_4)\psi'(a_3)
\end{align*}
}%
\noindent
It can be seen that in each equation the first product corresponds to the flat path and the remaining ones 
correspond to a concetenation of a pyramid and a flat path.
One can verify that we obtain the same expressions if we apply the Leibniz formula directly
to the c-free cumulants (with outer blocks associated with $\varphi$ 
and inner blocks associated with $\psi$). For instance, if $i_1=i_3$, we have
\begin{align*}
\varphi(a_1a_2a_3)&= R_2(a_1,a_3)r_1(a_2)+R_1(a_1)R_1(a_2)R_1(a_3)\\
&=(\varphi(a_1a_3)-\varphi(a_1)\varphi(a_2))\psi(a_2)+\varphi(a_1)\varphi(a_2)\varphi(a_3)
\end{align*}
where symbols $R_n$ and $r_n$ denote c-free and free cumulants, respectively. Thus,
\begin{align*}
\varphi(a_1a_2a_3)'&=(\varphi(a_1a_3)\varphi(a_1)-\varphi(a_1)\varphi(a_3))'\psi(a_2)+\\
&\quad +(\varphi(a_1a_2)-\varphi(a_1)\varphi(a_2))\psi'(a_2)\\
&\quad +\varphi(a_1)'\varphi(a_2)\varphi(a_3)+\varphi(a_1)\varphi'(a_2)\varphi(a_3)+\varphi(a_1)\varphi(a_2)\varphi'(a_3)\\
&=\varphi(a_1a_2)\psi'(a_2)+\varphi(a_1)'\varphi(a_2)\varphi(a_3),
\end{align*}
where we used $\varphi(a_1)=\psi(a_2)=\psi(a_3)=0$.
}
\end{Example}

We can also compute higher order derivatives of functionals $\widetilde{\Phi}_w$.
The combinatorics of the formula for the $m$-th derivative looks similar to that for the 
functionals $\Phi_{w}$, but there is one important difference: not all local maxima at level $1$ 
are treated in the same way since the first variable is $\varphi_{i_1}$-centered and the remaining ones are
$\psi_{i_k}$-centered, where $k>1$.

\begin{Theorem}
Let $a_1\in {\mathcal A}_{i_1}^{\circ},a_2\in {\mathcal A}_{i_2}^{\square}, \ldots, a_n\in \mathcal{A}_{i_n}^{\square}$, 
where $i_1\neq \cdots \neq i_n$ and let $w=j_1\cdots j_n\in \mathpzc{M}_{n}$.
Let $\pi(w)=\{V_1, \ldots, V_{r}\}$, where $V_1=\{k_1\}, \ldots, V_p=\{k_p\}$ are singleton blocks 
associated with colors $>1$ and perhaps $\{1\}$, if it is a singleton in $w$, and $V_{p+1}, \ldots, V_r$ are other blocks, 
where $p\leq r$.
\begin{enumerate}
\item  If $p\leq m$, then it holds that
\begin{align*}
\;\;\;\;\;\widetilde{\Phi}_{w}^{(m)}(a_1, \ldots, a_n)&=
\sum_{\stackrel{m_1+\cdots +m_r=m}{\scriptscriptstyle m_{k}\geq 1\;{\rm for}\;k\leq p}}
\frac{m!}{m_1! \cdots m_r!}
\prod_{s=1}^{p}
\eta_{i_{k_s}}^{(m_s)}(a_{k_s})
\prod_{s=p+1}^{r} \gamma_{|V_s|}^{(m_s)}((a_{k})_{k\in V_s})
\end{align*}
whenever $w$ is adapted to $(i_1, \ldots, i_n)$. 
\item
If $w$ is not adapted to $(i_1, \ldots, i_n)$ or $p>m$, then 
$\widetilde{\Phi}_{w}^{(m)}(a_1, \ldots, a_n)=0$.
\end{enumerate}
\end{Theorem}
\vspace{3pt}
{\it Proof.}
The proof is similar to that of Theorem 8.1. 
The first difference is that we have $\varphi_{i,t}$ at tensor positions of color $1$ and $\psi_{i,t}$ at 
tensor positions of colors $>1$. 
This difference is taken into account in the more general notation involving functionals 
$\eta_{i_k}$ and $\gamma_{|V_s|}$. The second difference is slightly more important: not all 
singletons (not all local maxima) are treated in the same way. 
If we apply the Leibniz formula to moments of order $n>1$ of $\widetilde{\Phi}_{w,t}$ 
given in the proof of Proposition 9.3, the role of functions $g_{k}(x)$ for $1\leq k \leq p$ 
used in the proof of Theorem 8.1 is played by 
$\varphi_{i_1,t}(a_{1})$ (if $\{1\}$ is a singleton of color $j_1=1$ determined by $w$) 
or $\psi_{i_k,t}(a_{k})$ (if $\{k\}$ is a singleton of color $j_k>1$) since 
in these cases 
{
\setlength{\belowdisplayshortskip}{\baselineskip}
\[
\varphi_{i_1,0}(a_1)=0\;\;\;{\rm and} \;\;\;\psi_{i_k}(a_k)=0.
\]
}%
\noindent
However, if $\{k\}$ is another singleton of color $j_k=1$, thus $k>1$ and, in general, we cannot 
claim that 
{
\setlength{\belowdisplayshortskip}{\baselineskip}
\[
\varphi_{i_k}(a_{k})=0 \;\;{\rm since}\;\;a_k\in {\rm Ker}(\psi_{i_k}).
\]
}%
\noindent
In consequence, the role of $g_{k}(x)$ for $k>p$ is played by 
$\varphi_{i_k}(a_{k})=\beta_{1}(a_{k})$ for $k>1$ whenever $\{k\}$ is a singleton, and 
by the remaining Boolean cumulants $\beta_{|V_s|}((a_{k})_{k\in V}$ 
arising from the excursions. 
In other words, singletons at level $j>1$ and possibly only one singleton at level $1$ 
(namely, $\{1\}$) are counted as singletons in the way we partition the number $m$ 
since we have to assign a derivative to these singletons, whereas the remaining singletons at level $1$ 
do not need a derivative since their moments of order $1$ are computed with respect to 
$\varphi_{i_k}$ and they may not be equal zero.
Of course, we need to take into account only those words $w\in \mathpzc{M}_{n}$ 
which are adapted to $\ell=(i_1, \ldots, i_n)$ since only for these words the moment
functional $\widetilde{\Phi}_{w,t}(a_1, \ldots , a_n)$ may not be equal to zero. This completes the proof. 
\hfill $\blacksquare$\\


\begin{thebibliography}{99}
\bibitem{Av} 
D.~Avitzour, Free products of $C^{*}$- algebras,
{\it Trans.~Amer.~Math.~Soc.} {\bf 271} (1982), 423-465.
\bibitem{BSh}
S.~T.~Belinschi, D.~Shlyakhtenko, Free probability of type B: analytic interpretation and applications, {\it Amer. J. Math.} {\bf 134} 
(2012) no. 1, 193–234.
\bibitem{BGN}
Ph.~Biane, F.~Goodman, A.~Nica, Alexandru Non-crossing cumulants of type B. {\it Trans. Amer. Math. Soc.} 355 (2003), no. 6, 2263–2303. 
\bibitem{BLS}
M.~Bo\.{z}ejko, M.~Leinert, R.~Speicher, Convolution and limit theorems for conditionally free random variables, 
{\it Pacific J. Math.} {\bf 175} (1996), no. 2, 357–388.
\bibitem{CDG}
G.~C\'ebron, A.~Dahlqvist, F.~Gabriel, Freeness of type B and conditional freeness for random matrices. {\it Indiana Univ. Math. J.} 73 (2024), no. 3, 1207–1252.
\bibitem{FMNS}
M.~F\'evrier, M.~Mastnak, A.~Nica, K.~Szpojankowski, A construction which relates c-freeness to infinitesimal freeness, 
{\it Adv. Appl. Math.} {\bf 110} (2019), 299–341.
\bibitem{FN}
M.~F\'evrier, A.~Nica, Infinitesimal non-crossing cumulants and free probability of type B, {\it J. Funct. Anal.} {\bf 258} 
(2010), no. 9, 2983–3023.
\bibitem{F}
M.~F\'evrier, Higher order infinitesimal freeness. {\it Indiana Univ. Math. J.} {\bf 61} (2012), no. 1, 249–295.
\bibitem{FH}
K.~Fujie, T.~Hasebe, Free probability of type B prime, {\it Trans. Amer. Math. Soc.} {\bf 378} (2025), no. 8, 5551–5577.
\bibitem{L1}
R.~Lenczewski, Unification of independence in quantum probability,
{\it Infin.~Dimens. Anal.~Quantum Probab. Relat.~Top.} {\bf 1} (1998), 383-405.
\bibitem{L2} 
R.~Lenczewski, Reduction of free independence to tensor independence,
{\it Infin. Dimens. Anal. Quantum Probab. Relat. Top.} {\bf 7} (2004), 337-360.
\bibitem{L3}
R.~Lenczewski, Decompositions of the free additive convolution, 
{\it J. Funct. Anal.} {\bf 246} (2007), 330-365.
\bibitem{L4}
R.~Lenczewski, Motzkin path decomposition of functionals in noncommutative probability, 
{\it Infin.~Dimens. Anal.~Quantum Probab. Relat.~Top.}, {\bf 25} (2022), no. 4, Paper No. 2240002.
\bibitem{L5}
R.~Lenczewski, Decomposition of free cumulants, ArXiv:2307:02281 (2023).
\bibitem{MT}
J.~Mingo, P.-L.~Tseng, Infinitesimal operators and the distribution of anticommutators and commutators, 
{\it J. Funct. Anal.} {\bf 287} (2024), no. 9, Paper No. 110591, 35 pp.
\bibitem{Sh}
D.~Shlyakhtenko, Free probability of type B and asymptotics of finite-rank perturbations of random matrices,
{\it Indiana Univ. Math. J.} {\bf 67} (2018), no. 2, 971–991.
\bibitem{T}
P.-L.~Tseng, A unified approach to infinitesimal freeness with amalgamation. {\it Internat. J. Math.} {\bf 34} 
(2023), no. 13, Paper No. 2350079, 26 pp. 
\bibitem{V1}
D.~Voiculescu, Symmetries of some reduced free product $C^{*}$-algebras, Operator Algebras and Their Connections with Topology and Ergodic Theory, Lecture Notes in Mathematics, Vol. 1132, Springer Verlag, 1985, pp. 556-588.
\bibitem{V2}
D.~Voiculescu, Limit laws for random matrices and free products, {\it Invent. Math.} {\bf 104} (1991), 201-220.

\end{thebibliography}
\end{document}